\documentclass[12pt,a4paper]{article}   

\usepackage[twoside]{geometry}

\usepackage{calc}
\usepackage[latin1]{inputenc}
\usepackage[T1]{fontenc}
\usepackage[all,12pt]{xy}

\usepackage[frenchb,english]{babel}

\usepackage{mdwtab} 
\usepackage{delarray} 
\usepackage{hhline}
\usepackage{longtable}
\usepackage{multicol}
\usepackage{xspace}

\usepackage{12many}
\newOTMstyle{dblbracket}{[\hspace{-2pt}[ #1,#2 ]\hspace{-2pt}]} 
\newOTMstyle{CeilFloor}{\lceil #1,#2 \rfloor}
\setOTMstyle{CeilFloor}

\usepackage{amsmath}
\usepackage{amssymb}

\newcommand{\Cal}{\mathcal}

\usepackage{accents}

\newcommand{\BB}[1]{\mathbb{#1}}

\newcommand{\FF}{\BB{F}}
\newcommand{\RR}{\BB{R}}
\newcommand{\NN}{\BB{N}}

\newcommand{\cube}{\square}
\newcommand\he {{\,\simeq\,}}

\DeclareMathOperator{\hocolim}{hocolim}

\DeclareMathOperator{\colim}{colim}

\DeclareMathOperator{\Mor}{Mor}
\DeclareMathOperator{\Ob}{Ob}
\DeclareMathOperator{\Top}{Top}

\newcommand\In[1]{\accentset{\circ}{#1}}

\usepackage{pst-all}

\usepackage[thref,thmmarks]{ntheorem}

\theoremstyle{changebreak}

\newtheorem{theo}{Theorem}[section]
\newtheorem{THEO}[theo]{Theorem}
\newtheorem{prop}[theo]{Proposition}
\newtheorem{PROP}[theo]{Proposition}

\newtheorem{LEMM}[theo]{Lemma}

\newtheorem{CORO}[theo]{Corollary}

\theoremstyle{nonumberbreak}
\newtheorem{THEOA}{Main Theorem}

\theoremstyle{changebreak}
\theorembodyfont{\upshape}

\newtheorem{defi}[theo]{Definition}
\newtheorem{DEFI}[theo]{Definition}
\newtheorem{exem}[theo]{Example}
\newtheorem{EXEM}[theo]{Example}
\newtheorem{exems}[theo]{Examples}
\newtheorem{EXEMS}[theo]{Examples}

\theoremstyle{nonumberbreak}

\theoremstyle{nonumberplain}
\theoremheaderfont{\bfseries}\theorembodyfont{\upshape}
\theoremseparator{.}
\theoremsymbol{\hskip 2ex \ensuremath{\square}}
\newtheorem{demo}{Proof}[subsection]
\newtheorem{DEMO}{Proof}[subsection]

\usepackage{fancyhdr}
\begin{document}
\title{
Configuration spaces of an embedding torus\\ and cubical spaces\\
}
\author{Jean-Philippe~Jourdan}
\date{March 20, 2006}
\maketitle
\begin{abstract}
For a smooth manifold $M$ obtained as an embedding torus, $A \bigcup C\times[-1,1]$, we consider the ordered configuration space $\FF_k(M)$ of $k$ distinct points in $M$.
We show that there is a homotopical cubical resolution of $\FF_k(M)$ defined from the configuration spaces of $A$ and $C$.
From it, we deduce a universal method for the computation of the pure braid groups of a manifold. 
We illustrate the method in the case of the Möbius band.
\end{abstract}

\section*{Introduction}

For a manifold $M$, the configuration space of $k$ points in $M$ is the space
\[\FF_k(M) = \left\{ (x_1,\dots,x_k)\in M^{\times k} \mid i\neq j\Rightarrow x_i\neq x_j\right\}.\]
As pointed-out by R.~Fox and L.~Neuwirth \cite{Fox1962}, configuration spaces are closely related to braid groups. In fact the fundamental group of the configuration spaces of the euclidean plane $\RR^2$ is Artin's pure braid groups \cite{Artin1947}. Hence, the pure braid group on $k$ strands of a manifold $M$ is naturally defined as the fundamental group of $\FF_k(M)$.

The aim of this paper is two-folds.
On one hand, we give an explicit homotopy model for the configuration spaces of a smooth manifold $M$ obtained as an embedding torus described below, $A \bigcup C\times[-1,1]$.
The model is given in terms of a homotopy cubical space, or $h\cube$-space.
On the other hand, using the previous model, we give a method for the determination of the braid groups of a manifold.
To carry on such a calculation, one does not need to understand the full machinery behind the homotopy model. In fact, only basic knowledge on the fundamental group is required.

More precisely, 
let $A$ and $C$ be  two manifolds of dimension $n$ and $n-1$ respectively.
Denote by $\partial A$ the boundary of $A$ and fix two disjoint embeddings $i^-:C\to \partial A$ and $i^+:C\to \partial A$.
The homotopy colimit of the diagram
$\xymatrix{C\ar@<-.2em>[r]_{i^+}\ar@<.2em>[r]^{i^-}&A}$
is called an \emph{embedding torus} and is a topological $n$-manifold with corners:
\[\displaystyle M=A\bigcup_{i^-\times\{-1\}\bigsqcup i^+\times\{+1\}} C\times[-1,1]=\hocolim \left(\xymatrix{C\ar@<-.2em>[r]_{i^+}\ar@<.2em>[r]^{i^-}&A}\right).\]
In the sequel, the maps $i^-$ and $i^+$ will be denoted by
$i^\epsilon$, respectively with $\epsilon=-1$ and $\epsilon=+1$.
This situation is sufficiently general as the next example shows.
\begin{exems}\label{Def-exemples}
\begin{enumerate}
	\item  Let $C$  be the interval $[-1,1]$ and $A$ be the square $[-1,1]\times[-1,1]$.
\begin{itemize}
	\item  If the maps $i^\epsilon$ are defined by $i^\epsilon(x)=(x,\epsilon)$ then $M$ is a cylinder $\RR^2\setminus\{0\}$.
	\item  If the maps $i^\epsilon$ are defined by $i^\epsilon(x)=(\epsilon x,\epsilon)$ then $M$ is the Möbius band, denoted by $\Cal M$.
\end{itemize}
	\item The surfaces of genus $g\ge 1$, oriented and non-oriented, can be obtained with $A$ the disk $D^2$ with $2g$ small disks removed and $C$ the disjoint union of $g$ copies of the circle $S^1$.
\item If $U$ and $V$ are two $n$-manifolds, the connected sum  $U \# V$ can be obtained with $A=\left(U\setminus \In D^n\right) \bigsqcup \left(V\setminus \In D^n\right)$ and $C=S^{n-1}$.
	\item If $M$ is obtained from a manifold $N$ by a surgery, then $M$ can be obtained as an embedding torus with $A=(N\setminus S^r\times \In D^{n-r}) \bigsqcup (D^{r+1}\times S^{n-r-1})$ and $C=S^{r}\times S^{n-r-1}$.
\end{enumerate}
\end{exems}
\begin{THEOA}\label{TH-EXIST}
There exists a homotopical cubical space $X_\bullet^k$,
defined using the configuration spaces of $A$ and $C$ and the embeddings $i^-$ and $i^+$,
such that the geometrical realization $|X^k_\bullet|$ has the same homotopy type as $\FF_k(M)$.
\end{THEOA}
\begin{exem}
When $k=2$, the homotopical cubical space model of $\FF_2(M)$ is given by the diagram
\[
X^2_\bullet : 
\xymatrix{
\FF_2(C) \ar@<-0.5em>[rr]|-{d_2^-}\ar@<0.5em>[rr]|-{d_1^+}\ar@<-1.5em>[rr]|-{d_2^+}\ar@<1.5em>[rr]|-{d_1^-}	&&
C\times A \bigsqcup A\times C \ar@<-0.5em>[r]_-{d_1^+}\ar@<0.5ex>[r]^-{d_1^-}&
\FF_2(A)\\
}
\text{ with } d^{\epsilon}_1\circ d^{\epsilon'}_2 \simeq d^{\epsilon'}_1 \circ d^{\epsilon}_1.
\]
\end{exem}

The paper is organized as follows.
Section \ref{sec:CubicalSpaces} introduces the notion of homotopical cubical space.
The cubical resolution associated to the embedding torus $M$ is described in Section \ref{sec:CubeResol}. Finally, in Section \ref{sec:Braids}, we explain the method of computation of the braid groups of $M$ and detail the case of the Möbius band.

All spaces considered here are assumed to be compactly generated and with the homotopy type of a CW-complex. We denote by $\Top$ the category of such spaces.

\medbreak
\noindent\textbf{Note: }
The section on braid groups can be read rather independently,
only Example \ref{realization-cube2} and SubSection \ref{Cube2resol} are required to understand it.

\section{Cubical Spaces}\label{sec:CubicalSpaces}

In this section, we define the notion of homotopical cubical spaces, or $h\cube$-spaces.
Those spaces are defined as homotopy $\cube$-diagrams for a well chosen category $\cube$.

\subsection{Homotopy limits and colimits}\label{Homotopy-Lim-Colim}

We first expose some results of R.~Vogt \cite{Vogt1977,Vogt1973} concerning homotopy colimit of up to homotopy  commutative  diagram of spaces.
 
Let  $\Cal C$ be a small category. 
\begin{DEFI}\label{def:hcdiagram}
For each $A\in\Cal C$ and $B\in\Cal C$, let
\[
\begin{array}{l}
\Cal C_n(A,B)=\{(f_n,f_{n-1},\dots,f_1)\in(\Mor \Cal C)^n \mid f_n\circ\cdots\circ f_1 : A\to B \text{ in } \Cal C\}, \quad n\in\NN^*,\cr
\Cal C_0(A,B)=\begin{array}\{{ll}. 
				\{id_A\}	&	A=B	\cr
				\emptyset	&	A\neq B.\cr
				\end{array}
\end{array}
\]
A \emph{homotopy $\Cal C$-diagram $D$}, or \emph{$h\Cal C$-diagram}, consists of 
\begin{itemize}
	\item a map $D_0 : \Ob\Cal C\to \Top$,
	\item a collection of maps indexed by  $\Ob\Cal C$
	\[D_B:\bigsqcup_{n\geq 0}\bigsqcup_{A\in\Cal C} \Cal C_{n+1}(A,B)\times I^n\times D_0 A\to D_0B\]
	where $I=[0,1]$ and such that:
	\begin{enumerate}
	\item $D_B(id_B;x)=x$, for $x\in D_0B$.
	\item For $n>0$, $(f_{n},f_{n-1},\dots,f_0)\in\Cal C_{n+1}(A,B)$ and $(t_{n},t_{n-1},\dots,t_1)\in I^n$,
	\[
	\begin{array}{l}
	D_B(f_{n},t_{n},f_{n-1},t_{n-1},\dots,f_0;x)\cr
	\qquad=\begin{array}\{{lll}.
	D_B(f_{n},t_{n},\dots,f_1;x)	&	f_0=id,	\cr
	D_B(f_{n},t_{n},\dots,f_{i+1},t_{i+1}t_{i},f_{i-1},\dots,f_0;x)	& f_i=id, & 0<i<n,	\cr
	D_B(f_{n-1},t_{n-1},\dots,f_0;x)	&	f_n=id,	\cr
	D_B(f_{n},t_{n},\dots,f_{i+1},t_{i+1},f_{i}\circ f_{i-1},t_{i-1},\dots,f_0;x)	&t_i=1,\cr
	D_B(f_{n},t_{n},\dots,f_{i};D_C(f_{i-1},t_{i-1},\dots,f_0;x))	&	t_i=0,\cr
	\end{array}\end{array}
	\]
	with $C$ the source of $f_{i}$.
	\end{enumerate}
\end{itemize}
\end{DEFI}
When no confusion is possible, we write  $D$ in place of $D_B$.
Roughly, an $h\Cal C$-diagram is a kind of functor from $\Cal C$ to $\Top$ in which the equality between $D(f\circ g)$ and $D(f)\circ D(g)$ is valid only up to a coherent homotopy, see \cite[Example page~18]{Vogt1973}.
The definition of the homotopy colimit of such a diagram states as follows.
\begin{DEFI}\label{Def-colimitehomotopique}
Let $D$ be a $h\Cal C$-diagram. The \emph{homotopy colimit} of $D$, denoted $\hocolim D$, is defined as the space
\[
\hocolim D= \bigsqcup_{A,B\in\Cal C}\bigsqcup_{n\geq 0} \Cal C_n(A,B)\times I^n\times D_0 A \Big/ \sim
\]
where the relation $\sim$ is given by
	\[
	\begin{array}{l}
	(t_{n},f_{n},\dots,t_1,f_1;x)\cr
	\qquad=\begin{array}\{{lll}.
	(t_{n},f_{n},\dots,t_2,f_2;x)	&	f_1=id,	\cr
	(t_{n},f_{n},\dots,f_{i+1},t_{i}t_{i-1},f_{i-1},\dots,f_1;x)	& f_i=id, & 1<i,	\cr
	(t_{n},f_{n},\dots,t_{i+1},f_{i+1}\circ f_{i},t_{i-1},\dots,f_1;x)	&t_i=1, & i<n,\cr
	(t_{n-1},f_{n-1},\dots,f_1;x)	&t_n=1,\cr
	(t_{n},f_{n},\dots,f_{i+1};D_C(f_{i},t_{i-1},\dots,f_1;x))	&	t_i=0,\cr
	\end{array}\end{array}
	\]
	with $C$ the source of $f_{i}$.\\
For $p\in\NN$, the image of $\displaystyle\bigsqcup_{A,B\in\Cal C}\bigsqcup_{n=0}^p \Cal C_n(A,B)\times I^n\times D_0 A$ in $\hocolim D$ is denoted by $\hocolim^p D$.
\end{DEFI}

\subsection{The $\cube$ category}\label{CubeCat}

In this paragraph, we define the category $\cube$.
We first introduce the sets $\Cal C^n_p$, used in the definition of the morphisms of $\cube$.
\begin{DEFI}
Let $\ito{n}$ be the set of integers between $1$ and $n$ and
$$\Cal C^n_p =\left\{ P\subset \ito{n} \mid  |P| = p \right\}=\left\{ (i_1,\dots,i_p)\in \ito{n}^p \mid  i_1<\dots<i_p \right\}$$
be the set of subsets of $\ito{n}$ of cardinal $p$.
We define the map $\complement : \Cal C^n_p	\to \Cal C^n_{n-p}$ that sends an element $i=(i_1,\dots,i_p)$ to its complement in $\ito{n}$.
We also define two associative operations $\wedge$ and $\vee$:\\
For  $p<q<n$,\quad $\wedge :
 \begin{array}[t]{rcl}
\Cal C^p_q \times \Cal C^n_p	&\to& 	\Cal C^n_q	\cr
( j, i)=((j_1,\dots,j_q),(i_1,\dots,i_p))	&\mapsto&	  j \wedge  i=(i_{j_1},\dots,i_{j_q})
\end{array}$\\[.5em]
and for $1\leq p+q \leq k$,\quad
$\vee :
\begin{array}[t]{rcl}
\Cal C^{n-p}_{q} \times \Cal C^n_p	&\to& 	\Cal C^n_{p+q}	\cr
( j, i)	&\mapsto& j \vee  i	=  i + ( j \wedge \complement i)
\end{array}
$\\
where $+$ stands for the union of sets in $\ito{n}$.
At last, we define a map $[\quad,\quad]$ acting like a partial inverse for the $\vee$ operation:
\[[\quad,\quad] :
\begin{array}[t]{rcl}
\Cal C^{n-p}_{q} \times \Cal C^n_p	&\to& 	\Cal C^{p+q}_{p}	\cr
( j, i)	&\mapsto& [ j , i] 
\end{array}
\]
where $[ j , i]$ is the only element such that $[ j , i]\wedge ( j \vee  i)=  i$.
\end{DEFI}

\begin{EXEMS}
\begin{itemize}
	\item If $i=(2,\textbf{3},5,\textbf{7},9,\textbf{11},13,17)\in\Cal C_8^{18}$ and $j=(2,4,6)\in\Cal C_3^8$, then\\*
	\centerline{$j\wedge i=(i_2,i_4,i_6)=(3,7,11)\in\Cal C_3^{18}$.}
	\item If $i=(2,3,5,7,9,11,13,17)\in\Cal C_8^{18}$ and $j=(1,4,6,9)\in\Cal C_4^{10}$, then\\*
			\centerline{$\complement i=(1,4,6,8,10,12,14,15,16,18)$ thus $j \wedge \complement i=(1,8,12,16)$.}
			The set $j \vee  i$ is given by
			$\begin{array}[t]{rcl}
			i + ( j \wedge \complement i)&=&(2,3,5,7,9,11,13,17)+(1,8,12,16)\cr
					&=&(1,2,3,5,7,8,9,11,12,13,16,17)\in\Cal C_{12}^{18}.\cr\end{array} $\\
			In order to precise the set $[j,i]\in\Cal C_8^{12}$, we have to determine the position of the integers forming  $i=(2,3,5,7,9,11,13,17)$ in $j \vee  i=(1,\textbf{2},\textbf{3},\textbf{5},\textbf{7},8,\textbf{9},\textbf{11},12,\textbf{13},16,\textbf{17})$. We deduce $[j,i]=(2,3,4,5,7,8,10,12)$.
\end{itemize}
\end{EXEMS}
We can use the operation $[\quad,\quad]$ to extend the definition of $\vee$ to a coloring of the integers forming $i$ and $j$.
\begin{DEFI}
For $ i=(i_1,\dots,i_p)\in\Cal C^n_p$, $\epsilon=(\epsilon_1,\dots,\epsilon_p)\in\{-1,1\}^p$, $ j=(j_1,\dots,j_q)\in\Cal C^{n-p}_q$ and $\epsilon'=(\epsilon'_1,\dots,\epsilon'_q)\in\{-1,1\}^q$, we define $\epsilon''=(\epsilon''_1,\dots,\epsilon''_{p+q})\in\{-1,1\}^{p+q}$
as
\[\begin{array}{lll}
\epsilon''_{l_\alpha}=\epsilon_\alpha & \text{with }l=[j,i]\in\Cal C_p^{p+q} &\text{and } \alpha\in\ito{p},\cr
\epsilon''_{s_\beta}=\epsilon'_\beta & \text{with } s=\complement l\in\Cal C_q^{p+q}	&\text{and }	\beta\in\ito{q}.\cr
\end{array}
\]
When no confusion is possible, the element $\epsilon''\in\{-1,1\}^{p+q}$ is denoted by $\epsilon_{ j\vee i}$.
\end{DEFI}
\begin{DEFI}
The category $\cube$ is defined by
\[\begin{array}[t]{l}
\Ob \cube = \NN	\cr
\Mor_\cube (p,r) =	
\begin{array}[t]\{{lll}.
	\left\{ f_{{i}}^{\epsilon_{ i}} \mid  i \in \Cal C_p^{p-r},	\epsilon_{ i}\in\{-1,1\}^{p-r}\right\}&	p>r, \text{ with } f_{{j}}^{\epsilon_{ j}}\circ f_{{i}}^{\epsilon_{ i}}=f_{ j \vee i}^{\epsilon_{ j\vee i}},\\
	\{id_p\}						&	p=r,\\
	\emptyset						&	p<r.	\\
\end{array}
\end{array}\]
For $p\in\NN$, we also define $\cube_p$ as the full subcategory of $\cube$ whose objects are $\oto{p}$.

\end{DEFI}
We finally prove some assumptions on the operations defined on the spaces $\Cal C^*_*$. These claims are convenient for proving the existence of a $h\cube$-space structure.
\begin{PROP}
\begin{itemize}
	\item [(1)] For $j\in\Cal C^{p}_{q}$ and $i\in  \Cal C^{n}_{p}$,
\quad $\complement( j\wedge  i) =(\complement j) \vee (\complement i)$.
	\item [(2)] For $j\in\Cal C^{n-p}_{q}$ and $i\in  \Cal C^{n}_{p}$,
\quad $ (\complement j) \wedge (\complement i) =\complement( j\vee  i)$.
	\item [(3)] For $j\in\Cal C^{n-p}_{q}$ and $i\in  \Cal C^{n}_{p}$, $[j,i]$ is the only element of $\Cal C^{p+q}_{p}$ which satisfies the relation
\[\complement [j,i]\wedge(j\vee i) = j\wedge \complement i .\]
	\item [(4)] For $k\in\Cal C^{n-p-q}_{r}$, $j\in\Cal C^{n-p}_{q}$ and $i\in  \Cal C^{n}_{p}$,
\[ [k\vee j,i]=[j,i]\wedge[k,j\vee i] \quad\text{and}\quad [k,j\vee i]=[k,j]\vee[k\vee j,i]. \]
\end{itemize}
\end{PROP}
\begin{DEMO}
Let $j\in\Cal C^{p}_{q}$ and $i\in  \Cal C^{n}_{p}$, the relations $(1)$ and $(2)$ follow from the equalities
\[ j\wedge  i = (j\wedge i)\cap i =  \complement(\complement j\wedge  i) \cap i
=\complement(\complement j\wedge  i + \complement i)
=\complement(\complement j\vee\complement i).
\]

Let $k\in\Cal C^{n-p-q}_{r}$, $j\in\Cal C^{n-p}_{q}$ and $i\in  \Cal C^{n}_{p}$.
Observe that $j\wedge \complement i=(j\vee i)\setminus i$ and $\complement k\wedge(j\vee i) = (j\vee i)\setminus (k\wedge(j\vee i))$. Hence, we have the equivalence
\[k\wedge(j\vee i) = i 
\Longleftrightarrow \complement k\wedge(j\vee i) = (j\vee i)\setminus i = j\wedge \complement i.\]
This shows the relation $(3)$.

Let $k\in\Cal C^{n-p-q}_{r}$, $j\in\Cal C^{n-p}_{q}$ and $i\in  \Cal C^{n}_{p}$, we have
\[
[j,i]\wedge[k,j\vee i]\wedge(k\vee j\vee i)=[j,i]\wedge (j\vee i) = i.
\]
The first equality of $(4)$ follows from the unicity of the element $[k\vee j,i]$. Furthermore, we have
\[
\complement[k,j]\wedge\complement [k\vee j,i]\wedge (k\vee j \vee i)
=\complement[k,j] \wedge(k\vee j)\wedge \complement i
=k\wedge\complement j \wedge \complement i = k\wedge\complement( j \vee  i).
\]
Thus, the second  equality of $(4)$ follows from  $(3)$.
\end{DEMO}

\subsection{$h\cube$-spaces}\label{CubeSpaces}

\begin{DEFI}
A \emph{homotopical cubical space}, or \emph{$h\cube$-space}, is defined as a  $h\cube$-diagram.
If $X_\bullet$ is a $h\cube$-space, we depict this diagram as a cubical complex
$$\xymatrix{
\ar@<.5em>[r]^<{d_1^-}\ar@<-.5em>[r]_<{d_5^+}\ar@3{{}{.}{}}[r]& X_{4}\ar@<.5em>[r]^{d_1^-}\ar@<-.5em>[r]_{d_4^+}\ar@3{{}{.}{}}[r]&
X_{3}\ar@<.5em>[r]^{d_1^-}\ar@<-.5em>[r]_{d_3^+}\ar@3{{}{.}{}}[r]&	 X_{2}\ar@<.5em>[r]^{d_1^-}\ar@<-.5em>[r]_{d_2^+}\ar@<-.17em>[r]\ar@<.17em>[r]&
X_{1}\ar@<.2em>[r]^-{d_1^-}\ar@<-.2em>[r]_-{d_1^+}&
X_{0}\\
}$$
with the conventions $X_i=X(i)$ and $d_i^\epsilon=X(f_i^\epsilon):X_p\to X_{p-1}$ for  $f_i^\epsilon\in\Mor(p,p-1)$. Finally \emph{the geometrical realization} of $X_\bullet$, is defined as $\hocolim X_\bullet$ and denoted by $|X_\bullet|$.
We define analogously a $h\cube_p$-space as a  $h\cube_p$-diagram.
We also denote by $X_{\leq p}$ the image of $\cube_p$ by $X$.
Observe that the new diagram $X_{\leq p}$ inherits a $h\cube_p$-space structure from the $h\cube$-space $X_\bullet$.
\end{DEFI}
The geometrical realization of $X_\bullet$ can be obtained by induction in the following way (see Segal \cite{Segal1974} for the simplicial case).
\begin{prop}[Inductive construction of the geometrical realization]\label{cubiq-induct}
Let $X_\bullet$ be $h\cube$-space.
For every $i\in\NN^*$, there exists a map
$\Phi_i: \partial([-1,1]^i)\times X_i\to |X_{\leq i}|$ such that
the square 
\[\xymatrix{
\partial([-1,1]^i)\times X_i\ar[r]^-{\Phi_i}\ar@{^(->}[d]	&	|X_{\leq i-1}|\ar@{^(->}[d]	\\
[-1,1]^i\times X_i\ar[r]	&	|X_{\leq i}|		\\
}\]
is a pushout.
Furthermore, the geometrical realization of $X_\bullet$ is the colimit of
\[
|X_{\leq 1}| \to |X_{\leq 2}| \to \cdots \to |X_{\leq i}| \to |X_{\leq i+1}| \to \cdots
\]
\end{prop}
\begin{EXEM}[Geometrical realization of a $h\cube_2$-space]\label{realization-cube2}
Start with a $h\cube_2$-space
\[X_\bullet : \xymatrix{
X_2 \ar@<-1em>[r]_{d_{2}^+}\ar@<-.4em>[r]|{d_{2}^-}\ar@<+.4em>[r]|{d_{1}^+}\ar@<1em>[r]^{d_{1}^-}
 & X_1 \ar@<-.3em>[r]_{d_{1}^+}\ar@<.3em>[r]^{d_{1}^-}	&	X_0\\
}.\]
From the definition of homotopy colimit (\ref{Def-colimitehomotopique}), we have $|X_{\leq0}|=X_0$ and
\[\begin{array}{rcl}
|X_{\leq 1}|&=&
\hocolim\left(\xymatrix{ X_1 \ar@<-.3em>[r]_{d_{1}^+}\ar@<.3em>[r]^{d_{1}^-}	&	X_0 }\right)\\
 &=& \displaystyle X_1\times[-1,1] \bigcup_{(x_1,\epsilon)\sim d_1^\epsilon(x_1)} X_0\\
&=&\colim\left(X_1 \leftarrow X_1\times\{-1,1\} \rightarrow |X_{\leq0}|\right).\end{array}\]

The $h\cube_2$-space $X_\bullet$ is the image of the category $\cube_2$ represented below
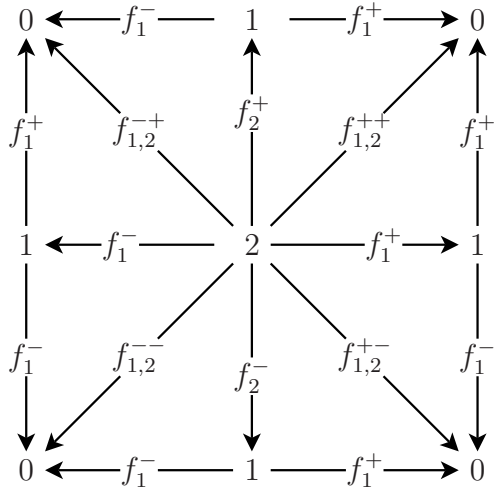
\begin{center}
\psset{xunit=1mm,yunit=1mm,runit=1mm}
\psset{linewidth=0.3,dotsep=1,hatchwidth=0.3,hatchsep=1.5,shadowsize=1}
\psset{dotsize=0.7 2.5,dotscale=1 1,fillcolor=black}
\psset{arrowsize=1 2,arrowlength=1,arrowinset=0.25,tbarsize=0.7 5,bracketlength=0.15,rbracketlength=0.15}
\begin{pspicture}(0,0)(65,65)
\rput(32.5,50){$f_2^+$}
\rput(32.5,15){$f_2^-$}
\rput(15,32.5){$f_1^-$}
\rput(50,32.5){$f_1^+$}
\rput(62.5,47.5){$f_1^+$}
\rput(47.5,62.5){$f_1^+$}
\rput(47.5,2.5){$f_1^+$}
\rput(2.5,47.5){$f_1^+$}
\rput(17.5,2.5){$f_1^-$}
\rput(17.5,62.5){$f_1^-$}
\rput(2.5,17.5){$f_1^-$}
\rput(62.5,17.5){$f_1^-$}
\rput(47.5,47.5){$f_{1,2}^{++}$}
\rput(47.5,17.5){$f_{1,2}^{+-}$}
\rput(17.5,17.5){$f_{1,2}^{--}$}
\rput(17.5,47.5){$f_{1,2}^{-+}$}
\psline[arrowscale=1.5 1.5]{->}(32.5,52.5)(32.5,60)
\psline[arrowscale=1.5 1.5]{->}(15,62.5)(5,62.5)
\psline[arrowscale=1.5 1.5]{<-}(2.5,60)(2.5,50)
\psline[arrowscale=1.5 1.5]{<-}(5,32.5)(12.5,32.5)
\psline[arrowscale=1.5 1.5]{->}(50,62.5)(60,62.5)
\psline[arrowscale=1.5 1.5]{<-}(62.5,60)(62.5,50)
\psline[arrowscale=1.5 1.5]{<-}(60,32.5)(52.5,32.5)
\psline[arrowscale=1.5 1.5]{->}(50,50)(60,60)
\psline[arrowscale=1.5 1.5]{<-}(5,60)(15,50)
\psline[arrowscale=1.5 1.5]{->}(50,15)(60,5)
\psline[arrowscale=1.5 1.5]{<-}(62.5,5)(62.5,15)
\psline[arrowscale=1.5 1.5]{<-}(60,2.5)(50,2.5)
\psline[arrowscale=1.5 1.5]{->}(32.5,12.5)(32.5,5)
\psline[arrowscale=1.5 1.5]{->}(15,2.5)(5,2.5)
\psline[arrowscale=1.5 1.5]{<-}(2.5,5)(2.5,15)
\psline[arrowscale=1.5 1.5]{->}(15,15)(5,5)
\rput(32.5,32.5){$2$}
\rput(62.5,32.5){$1$}
\rput(32.5,62.5){$1$}
\rput(2.5,32.5){$1$}
\rput(32.5,2.5){$1$}
\rput(2.5,62.5){$0$}
\rput(2.5,2.5){$0$}
\rput(62.5,2.5){$0$}
\rput(62.5,62.5){$0$}
\psline(35,35)(45,45)
\psline(32.5,47.5)(32.5,35)
\psline(20,45)(30,35)
\psline(27.5,32.5)(17.5,32.5)
\psline(2.5,35)(2.5,45)
\psline(2.5,30)(2.5,20)
\psline(30,2.5)(20,2.5)
\psline(35,2.5)(45,2.5)
\psline(32.5,30)(32.5,17.5)
\psline(35,32.5)(47.5,32.5)
\psline(62.5,35)(62.5,45)
\psline(30,30)(20,20)
\psline(35,30)(45,20)
\psline(62.5,30)(62.5,20)
\psline(27.5,62.5)(20,62.5)
\psline(37.5,62.5)(45,62.5)
\end{pspicture}

\end{center}
where all of the four external edges of the form $\cube_1: \xymatrix{ 0  &f_1^-  &1 \ar'[r][rr] \ar'[l][ll] &f_1^+ &0 }$ are identified.
Observe that the structure of $h\cube_2$-space is equivalent to the existence of homotopies $d_1^- d_2^-\he d_1^-d_1^-$, $d_1^+ d_1^- \he d_1^-d_2^+$, $d_1^+ d_2^+ \he d_1^+d_1^+$ and $d_1^- d_1^+ \he d_2^-d_1^+$.
Hence for every $x_2\in X_2$, we can define a path $\phi_1$ of $|X_{\leq 1}|$ as the restriction at $x_2$ of the homotopy between $d_1^- d_2^-$ and $d_1^-d_1^-$. Similarly, define $\phi_2$, $\phi_3$ and $\phi_4$ as the restrictions of the homotopies $d_1^+ d_1^- \he d_1^-d_2^+$,
$d_1^+ d_2^+ \he d_1^+d_1^+$ and $d_1^- d_1^+ \he d_2^-d_1^+$.
Also denote by $id$ the identity on $[-1,1]$. One construct a loop $\Phi_2(x_2,\bullet)$ in $|X_{\leq 1}|$ as follows:
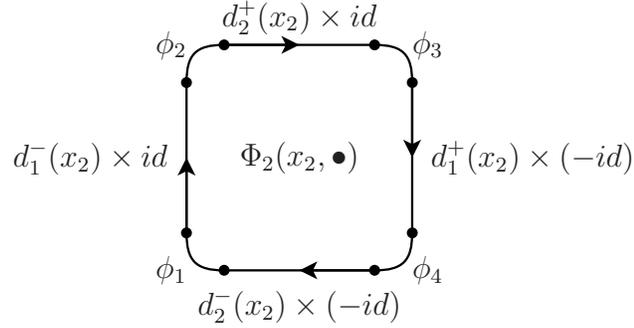
\begin{center}
\psset{xunit=1mm,yunit=1mm,runit=1mm}
\psset{linewidth=0.3,dotsep=1,hatchwidth=0.3,hatchsep=1.5,shadowsize=1}
\psset{dotsize=0.7 2.5,dotscale=1 1,fillcolor=black}
\psset{arrowsize=1 2,arrowlength=1,arrowinset=0.25,tbarsize=0.7 5,bracketlength=0.15,rbracketlength=0.15}
\begin{pspicture}(0,0)(37.5,37.5)
\rput[r](2.5,20){$d_1^-(x_2)\times id$}
\psline(10,35)(30,35)
\psdots[](10,35)
(30,35)
\psline(35,30)(35,10)
\psdots[](35,30)
(35,10)
\psline[arrowscale=1.5 1.5]{<-}(35,20)(35,30)
\psline(30,5)(10,5)
\psdots[](30,5)
(10,5)
\psline[arrowscale=1.5 1.5]{<-}(20,5)(30,5)
\psline[arrowscale=1.5 1.5]{->}(5,10)(5,20)
\psline(5,10)(5,30)
\psdots[](5,10)
(5,30)
\psline[arrowscale=1.5 1.5]{->}(10,35)(20,35)
\rput[t](20,2.5){$d_2^-(x_2)\times (-id)$}
\rput[B](20,37.5){$d_2^+(x_2)\times id$}
\rput[l](37.5,20){$d_1^+(x_2)\times (-id)$}
\rput[r](5,35){$\phi_2$}
\rput[r](5,5){$\phi_1$}
\rput[l](35,35){$\phi_3$}
\rput[l](35,5){$\phi_4$}
\psbezier(5,30)(5,33.75)(6.25,35)(10,35)
\psbezier(30,35)(33.75,35)(35,33.75)(35,30)
\psbezier(35,10)(35,6.25)(33.75,5)(30,5)
\psbezier(10,5)(6.25,5)(5,6.25)(5,10)
\rput(20,20){$\Phi_2(x_2,\bullet)$}
\end{pspicture}

\end{center}
More precisely, $\Phi_2(x_2,\bullet)$ is the following composition of paths:
\[\Phi_2(x_2,\bullet)=  \phi_1 * (d_1^-(x_2)\times id) * \phi_2 * (d_2^+(x_2)\times id) * \phi_3 * (d_1^+(x_2)\times (-id)) * \phi_4 * (d_2^-(x_2)\times (-id)).\]
We extend it as a map $\Phi_2:X_2\times\partial[-1,1]^2 \to |X_{\leq 1}|$. Then we have a description of $|X_{\leq 2}|$ as a pushout:
\[\xymatrix{
X_2\times\partial[-1,1]^2\ar[r]^-{\Phi_2}\ar[d] & |X_{\leq 1}|\ar[d]\\
X_2\times [-1,1]^2\ar[r] & |X_{\leq 2}|\\
}\]
\end{EXEM}
\begin{demo}[\ref{cubiq-induct}]Consider a $h\square_i$-space $X_\bullet$.
We describe the geometrical realization of $|X_{\leq i}|$ using the spaces $|X_{\leq i-1}|$ and $X_i$.
Define the space $\displaystyle A_i = \bigsqcup_{\substack{0\leq j \leq i\\  0\leq n \leq i}} \Cal C_n(i,j)\times I^n\times X_i \Big/ \sim$,
where the relations $\sim$ are given by
	\[
	\begin{array}{l}
	(t_{n},f_{n},\dots,t_1,f_1;x)\cr
	\qquad\sim \begin{array}\{{lll}.
	(t_{n},f_{n},\dots,t_2,f_2;x)	&	f_1=id,	\cr
	(t_{n},f_{n},\dots,f_{i+1},t_{i}t_{i-1},f_{i-1},\dots,f_1;x)	& f_i=id, & 1<i,	\cr
	(t_{n},f_{n},\dots,t_{i+1},f_{i+1}\circ f_{i},t_{i-1},\dots,f_1;x)	&t_i=1, & i<n,\cr
	(t_{n-1},f_{n-1},\dots,f_1;x)	&t_n=1.\cr
	\end{array}\end{array}
	\]
Observe that in each relation, the element $x\in X_i$ remains unmodified.
Since all the spaces are assumed to be in $\Top$, this implies
\[
A_i = \displaystyle\bigsqcup_{\substack{0\leq j \leq i\\  0\leq n \leq i}} \Cal C_n(i,j)\times I^n\times X_i \Big/\!\!\sim\ 
= \displaystyle\left( \bigsqcup_{\substack{0\leq j \leq i\\  0\leq n \leq i}} \Cal C_n(i,j)\times I^n \Big/\sim \right )\times X_i
= [-1,1]^i\times X_i.
\]
From Definition \ref{Def-colimitehomotopique} of the homotopy colimit,
we deduce that  $|X_{\leq i}|$ is obtained as the quotient space
$
A_i \bigsqcup |X_{\leq i-1}| \Big/ \sim
$
where the relation $ \sim $ is given by
\[(t_{n},f_{n},\dots,t_1,f_1;x)\sim (t_{n},f_{n},\dots,f_{i+1};D_C(f_{i},t_{i-1},\dots,f_1;x))\] 
whenever  $t_i=0$ with $C$ the source of $f_{i+1}$.
We can rephrase this result by saying that there exists a map $\Phi_i$ such that the square below is a homotopy pushout.
\[\xymatrix{
\partial([-1,1]^i)\times X_i\ar[r]^-{\Phi_i}\ar@{^(->}[d]	&	|X_{\leq i-1}|\ar@{^(->}[d]	\\
[-1,1]^i\times X_i\ar[r]	&	|X_{\leq i}|		\\
}\]
\end{demo}
As a direct consequence, we have the next result on the determination of the fundamental group of the geometrical realization of a $h\cube$-space.
\begin{PROP}\label{PI1-cubiq}
Let  $X_\bullet$ be a $h\cube$-space, then the following spaces have the same fundamental group:
\begin{itemize}
	\item $|X_\bullet|$,
	\item $|X_{\leq 2}|$,
	\item The geometrical realization of the square space
	$\xymatrix{X_2^{(0)} \ar@<-0.5ex>[r]\ar@<0.5ex>[r]\ar@<-1.5ex>[r]\ar@<1.5ex>[r]
	&X_1\ar@<-0.5ex>[r]\ar@<0.5ex>[r]&X_0\\}$.
\end{itemize}
with $X_2^{(0)}$ a $0$-skeleton of the space $X_2$.
\end{PROP}
A last property on $h\cube$-spaces is required in order to prove the Main Theorem.
\begin{prop} \label{DoubleCubiq}
Let $X$ be a $h\cube$-space of the form
\[\xymatrix{
\ar@<.5em>[r]^-{d_1^-}\ar@<-.5em>[r]_-{d_5^+}\ar@3{{}{.}{}}[r]&
\displaystyle\bigsqcup_{i+j=4} X_{i,j}\ar@<.5em>[r]^-{d_1^-}\ar@<-.5em>[r]_-{d_4^+}\ar@3{{}{.}{}}[r]&
\displaystyle\bigsqcup_{i+j=3} X_{i,j}\ar@<.5em>[r]^-{d_1^-}\ar@<-.5em>[r]_-{d_3^+}\ar@3{{}{.}{}}[r]&	\displaystyle\bigsqcup_{i+j=2} X_{i,j}\ar@<.5em>[r]^-{d_1^-}\ar@<-.5em>[r]_-{d_2^+}\ar@<-.17em>[r]\ar@<.17em>[r]&	\displaystyle X_{0,1}\bigsqcup X_{1,0}\ar@<.2em>[r]^-{d_1^-}\ar@<-.2em>[r]_-{d_1^+}&
X_{0,0}\\
}\]
such that each map $d_\alpha^\epsilon :  X_{p,q} \to X_{p,q-1}\bigsqcup X_{p-1,q}$ restricts to 
\[d_\alpha^\epsilon : X_{p,q} \to X_{p,q-1} \quad\text {if }\quad\alpha\leq q
	\quad\text{and}\quad d_\alpha^\epsilon : X_{p,q} \to X_{p-1,q}  \quad\text {if }\quad \alpha> q.
\]
Then, we can form the following diagram
\[(X'):\xymatrix{
\ar@<.5em>[rr]^{d_1^-}\ar@<-.5em>[rr]_{d_4^+}\ar@3{{}{.}{}}[rr]&&
X_{3,0}\ar@<.5em>[rr]^{d_1^-}\ar@<-.5em>[rr]_{d_3^+}\ar@3{{}{.}{}}[rr]&&	X_{2,0}\ar@<.5em>[rr]^{d_1^-}\ar@<-.5em>[rr]_{d_2^+}\ar@<-.17em>[rr]\ar@<.17em>[rr]&&	X_{1,0}\ar@<.2em>[rr]^{d_1^-}\ar@<-.2em>[rr]_{d_1^+}&&
X_{0,0}\\
\ar@<.5em>[rr]^{d_2^-}\ar@<-.5em>[rr]_{d_5^+}\ar@3{{}{.}{}}[rr]&&
X_{3,1}\ar@<.5em>[rr]^{d_2^-}\ar@<-.5em>[rr]_{d_4^+}\ar@3{{}{.}{}}[rr]\ar@<-.2em>[u]_{d_1^-}\ar@<.2em>[u]^{d_1^+}&&	X_{2,1}\ar@<.5em>[rr]^{d_2^-}\ar@<-.5em>[rr]_{d_3^+}\ar@<-.17em>[rr]\ar@<.17em>[rr]\ar@<-.2em>[u]_{d_1^-}\ar@<.2em>[u]^{d_1^+}&&
X_{1,1}\ar@<.2em>[rr]^{d_2^-}\ar@<-.2em>[rr]_{d_2^+}\ar@<-.2em>[u]_{d_1^-}\ar@<.2em>[u]^{d_1^+}&&
X_{0,1}\ar@<-.2em>[u]_{d_1^-}\ar@<.2em>[u]^{d_1^+}\\
\ar@<.5em>[rr]^{d_3^-}\ar@<-.5em>[rr]_{d_6^+}\ar@3{{}{.}{}}[rr]&&
X_{3,2}\ar@<.5em>[rr]^{d_3^-}\ar@<-.5em>[rr]_{d_5^+}\ar@3{{}{.}{}}[rr]\ar@<-.5em>[u]_{d_1^-}\ar@<.5em>[u]^{d_2^+}\ar@2{{}{.}{}}[u]&&	X_{2,2}\ar@<.5em>[rr]^{d_3^-}\ar@<-.5em>[rr]_{d_4^+}\ar@<-.17em>[rr]\ar@<.17em>[rr]\ar@<-.5em>[u]_{d_1^-}\ar@<.5em>[u]^{d_2^+}\ar@2{{}{.}{}}[u]&&
X_{1,2}\ar@<.2em>[rr]^{d_3^-}\ar@<-.2em>[rr]_{d_3^+}\ar@<-.5em>[u]_{d_1^-}\ar@<.5em>[u]^{d_2^+}\ar@2{{}{.}{}}[u]&&
X_{0,2}\ar@<-.5em>[u]_{d_1^-}\ar@<.5em>[u]^{d_2^+}\ar@2{{}{.}{}}[u]\\
&&
\ar@<-.5em>[u]_<{d_1^-}\ar@<.5em>[u]^<{d_3^+}\ar@3{{}{.}{}}[u]&&
\ar@<-.5em>[u]_<{d_1^-}\ar@<.5em>[u]^<{d_3^+}\ar@3{{}{.}{}}[u]&&
\ar@<-.5em>[u]_<{d_1^-}\ar@<.5em>[u]^<{d_3^+}\ar@3{{}{.}{}}[u]&&
\ar@<-.5em>[u]_<{d_1^-}\ar@<.5em>[u]^<{d_3^+}\ar@3{{}{.}{}}[u]&&\\
}\]
where each line and each column inherit a $h\cube$-space structure.
Let $|X_{j,\bullet}|$ (resp. $|X_{\bullet,j}|$) be the geometrical realization of the cubical space on the $j$-th column (resp. row).\\
Then, the geometrical realization of $X_\bullet$ is the same as the one of any of the  $h\cube$-spaces
\[\xymatrix@C=1.5em{
\ar@<.5em>[r]^-{d_1^-}\ar@<-.5em>[r]_-{d_4^+}\ar@3{{}{.}{}}[r]&
|X_{3,\bullet}|\ar@<.5em>[r]^-{d_1^-}\ar@<-.5em>[r]_-{d_3^+}\ar@3{{}{.}{}}[r]&	|X_{2,\bullet}|\ar@<.5em>[r]^-{d_1^-}\ar@<-.5em>[r]_-{d_2^+}\ar@<-.17em>[r]\ar@<.17em>[r]&	|X_{1,\bullet}|\ar@<.2em>[r]^-{d_1^-}\ar@<-.2em>[r]_-{d_1^+}&
|X_{0,\bullet}|\\
}
\ \text{ and }\ 
\xymatrix@C=1.5em{
\ar@<.5em>[r]^-{d_1^-}\ar@<-.5em>[r]_-{d_4^+}\ar@3{{}{.}{}}[r]&
|X_{\bullet,3}|\ar@<.5em>[r]^-{d_1^-}\ar@<-.5em>[r]_-{d_3^+}\ar@3{{}{.}{}}[r]&	|X_{\bullet,2}|\ar@<.5em>[r]^-{d_1^-}\ar@<-.5em>[r]_-{d_2^+}\ar@<-.17em>[r]\ar@<.17em>[r]&	|X_{\bullet,1}|\ar@<.2em>[r]^-{d_1^-}\ar@<-.2em>[r]_-{d_1^+}&
|X_{\bullet,0}|\\
}.\]
~\end{prop}
\begin{DEMO}   
The proof is a straightforward application of \cite[theorem 2.4]{Vogt1977} once we have observed that the diagram $(X')$ admits a structure of  $h\cube\times h\cube$-diagram (see \cite[Definition 2.1]{Vogt1977}) since it is defined from the $h\cube$-diagram $X$.
\end{DEMO}

\section{Cubical resolution associated to an embedding torus}\label{sec:CubeResol}

The aim of this section is to prove the Main Theorem.
In Section \ref{Cube2resol}, we give a geometric description of the resolution $X^k_\bullet$ whereas in Section \ref{CubeResol} we prove that the  $h\cube$-space $X^k_\bullet$ is well defined and augmented. The section \ref{DEMOTHA} is focused on the proof of the Main Theorem.

\subsection{$h\cube_2$-space associated to an embedding torus}\label{Cube2resol}

Recall that 
$$\Cal C^k_p =\left\{ P\subset \ito{k} \mid  |P| = p \right\}=\left\{ (i_1,\dots,i_p)\in \ito{k}^p \mid  i_1<\dots<i_p \right\}.$$
We set $X^k_p={\Cal C}^{k}_{p}\times{\FF}_{p}(C)\times {\FF}_{k-p}(A)$.
This space is naturally identified with the space of configurations of $k$ points 
in $A\bigsqcup C$ with exactly $p$ points in $C$.
In this section, we define some maps
$d_\alpha^\epsilon : X^k_p\to X^k_{p-1}$, $\alpha\in\ito{p}$, $\epsilon\in\{-1,+1\}$
and give a geometric interpretation of them. They will form a first step toward the definition of the resolution $X_\bullet$ of $\FF_k(M)$.

Recall that $A$ is a manifold with boundaries.
The collar neighborhood theorem \cite[Theorem 17.1]{Milnor1974} asserts the existence of an embedding $f:\partial A\times [0,2]\to A$ that extends the natural inclusion $\partial A\times \{0\}\hookrightarrow A$. Later on, we identify a point  $(a,t)\in\partial A\times [0,2]$ with its image by $f$.
This embedding allows us to define a map $r$ used in our construction.
Geometrically, this map crushes down  the collar neighborhood $\partial A\times [0,2]$ linearly onto $\partial A\times [1,2]$ leaving the complement unchanged.
\begin{prop}\label{retract}
There is an injective map $r : A\to A\setminus(\partial A\times [0,1[)$ which is a weak deformation retract of the natural inclusion.
\end{prop}
\begin{demo}
Define the map $R:A\times [0,2]\to A$ such that for $(a,s)\in A\times [0,2]$:
\[
\begin{array}{rll}
			&	R(a,s)=a 				&  \text{if } a\in A\setminus(\partial A\times [0,2[), \cr
 \text{and}	&	R((a,t),s)=f(a,s+(2-s)t/2) &	\text{if } (a,t)\in\partial A\times [0,2].\cr
\end{array}\]
The map $r=R(\bullet,1)$ fulfills the requirement.
\end{demo}

\begin{defi}\label{def-dalpha}
 For $1\leq \alpha \leq p$, the map $d_\alpha^\epsilon:  X^k_p\to X^k_{p-1}$ is defined by:
\[
\begin{array}{c}
\xymatrix{
{\Cal C}^{k}_{p}\times{\FF}_{p}(C)\times {\FF}_{k-p}(A) \ar[rr]	&& 
{\Cal C}^{k}_{p-1}\times{\FF}_{p-1}(C)\times {\FF}_{k-p+1}(A)\\
}\cr
\xymatrix{{\left((i_1,\dots,i_{p}),(c_1,\dots,c_{p}),(a_1,\dots,a_{k-p})\right)}\ar@{|->}[d]	\\&}\cr
\left((i_1,\ldotp\ldotp,\widehat{i_\alpha},\ldotp\ldotp,i_{p}), (c_1,\ldotp\ldotp,\widehat{c_\alpha},\ldotp\ldotp,c_{p}),
\left(r(a_1),\ldotp\ldotp,r(a_{i_\alpha-\alpha}),i^\epsilon(c_\alpha),r(a_{i_\alpha-\alpha+1}),\ldotp\ldotp,r(a_{k-p})\right)\right)\cr
\end{array}
\]
\end{defi}
We give now a geometric interpretation of the maps $d_\alpha^\epsilon: X^k_p\to X^k_{p-1}$.\\[.5em]
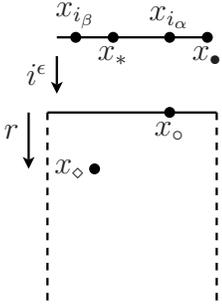
\begin{minipage}[c]{.2\textwidth}
\bigbreak
\psset{xunit=1mm,yunit=1mm,runit=1mm}
\psset{linewidth=0.3,dotsep=1,hatchwidth=0.3,hatchsep=1.5,shadowsize=1}
\psset{dotsize=0.7 2.5,dotscale=1 1,fillcolor=black}
\psset{arrowsize=1 2,arrowlength=1,arrowinset=0.25,tbarsize=0.7 5,bracketlength=0.15,rbracketlength=0.15}
\begin{pspicture}(0,0)(26.25,36.25)
\psline(3.75,25)(26.25,25)
\psline(5,35)(25,35)
\psdots[](7.5,35)
(20,35)
\psdots[](12.5,35)
(25,35)
\psdots[](10,17.5)
(20,25)
\psline[linestyle=dashed,dash=1 1](26.25,25)(26.25,0)
\psline[linestyle=dashed,dash=1 1](3.75,25)(3.75,0)
\rput[t](20,23.75){$x_\circ$}
\rput[r](8.75,17.5){$x_\diamond$}
\rput[b](7.5,36.25){$x_{i_\beta}$}
\rput[b](20,36.25){$x_{i_\alpha}$}
\rput[t](12.5,33.75){$x_\ast$}
\rput[t](25,33.75){$x_\bullet$}
\psline{->}(5,32.5)(5,27.5)
\psline{->}(1.25,25)(1.25,17.5)
\rput[r](3.75,30){$i^\epsilon$}
\rput[r](0,22.5){$r$}
\end{pspicture}

\end{minipage}
\begin{minipage}[c]{.79\textwidth}
Let $\left((i_1,\dots,i_{p}),(c_1,\dots,c_{p}),(a_1,\dots,a_{k-p})\right)\in X^k_p$.
As already stated, this object is identified with a configuration of $k$ points, $(x_1,\dots,x_k)$, living in $A\bigsqcup C$ with exactly $p$ of them in $C$.
 On the picture on the left, the map $i^\epsilon$ is given by the natural inclusion of $C$ (represented by the interval at the top of the picture) into $\partial A$ (represented by the other non dashed interval). The retraction $r$ sends down vertically the points living in $A$.
\end{minipage}\\[.5em]
The action of the map $d_\alpha^\epsilon$ on the above configuration consists of two successive steps:\\[.5em]
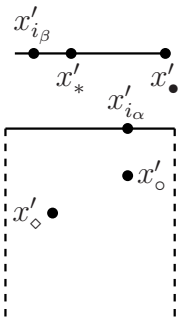
\begin{minipage}[c]{.2\textwidth}
\bigbreak
\psset{xunit=1mm,yunit=1mm,runit=1mm}
\psset{linewidth=0.3,dotsep=1,hatchwidth=0.3,hatchsep=1.5,shadowsize=1}
\psset{dotsize=0.7 2.5,dotscale=1 1,fillcolor=black}
\psset{arrowsize=1 2,arrowlength=1,arrowinset=0.25,tbarsize=0.7 5,bracketlength=0.15,rbracketlength=0.15}
\begin{pspicture}(0,0)(26.25,36.25)
\psline(3.75,25)(26.25,25)
\psline(5,35)(25,35)
\psdots[](7.5,35)
(20,25)
\psdots[](12.5,35)
(25,35)
\psdots[](10,13.75)
(20,18.75)
\psline[linestyle=dashed,dash=1 1](26.25,25)(26.25,0)
\psline[linestyle=dashed,dash=1 1](3.75,25)(3.75,0)
\rput[l](21.25,18.75){$x'_\circ$}
\rput[r](8.75,13.75){$x'_\diamond$}
\rput[b](7.5,36.25){$x'_{i_\beta}$}
\rput[b](20,26.25){$x'_{i_\alpha}$}
\rput[t](12.5,33.75){$x'_\ast$}
\rput[t](25,33.75){$x'_\bullet$}
\end{pspicture}

\end{minipage}
\begin{minipage}[c]{.79\textwidth}
\begin{enumerate}
	\item It pushes away from the boundaries all points in $A$ using the retraction~$r$,
	\item It sends the $\alpha$-th point living in $C$, $x_{i_\alpha}=c_\alpha$,  into $\partial A$ using the map~$i^\epsilon$.
\end{enumerate}
Let $(x'_1,\dots,x'_k)$ denotes the new configuration obtained,
in particular, $x'_{i_\alpha}=i^\epsilon(x_{i_\alpha})=i^\epsilon(c_\alpha)$.
Observe that the use of the retraction is necessary since one cannot ensure that the point $i^\epsilon(c_\alpha)$ is not already in the configuration $(x_1,\dots,x_k)$.
\end{minipage}\\[.5em]
Now if we assume $\beta>\alpha$, then $c_\beta=x'_{i_\beta}$ is the $(\beta-1)$-th point living in $C$ in the configuration $(x'_1,\dots,x'_k)$.\\*
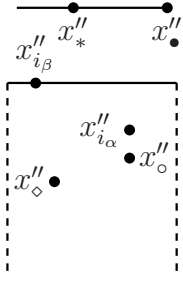
\begin{minipage}[c]{.2\textwidth}
\bigbreak
\psset{xunit=1mm,yunit=1mm,runit=1mm}
\psset{linewidth=0.3,dotsep=1,hatchwidth=0.3,hatchsep=1.5,shadowsize=1}
\psset{dotsize=0.7 2.5,dotscale=1 1,fillcolor=black}
\psset{arrowsize=1 2,arrowlength=1,arrowinset=0.25,tbarsize=0.7 5,bracketlength=0.15,rbracketlength=0.15}
\begin{pspicture}(0,0)(26.25,35)
\psline(3.75,25)(26.25,25)
\psline(5,35)(25,35)
\psdots[](7.5,25)
(20,18.75)
\psdots[](12.5,35)
(25,35)
\psdots[](10,11.88)
(20,15)
\psline[linestyle=dashed,dash=1 1](26.25,25)(26.25,0)
\psline[linestyle=dashed,dash=1 1](3.75,25)(3.75,0)
\rput[l](21.25,15){$x''_\circ$}
\rput[r](8.75,11.88){$x''_\diamond$}
\rput[b](7.5,26.25){$x''_{i_\beta}$}
\rput[r](18.75,18.75){$x''_{i_\alpha}$}
\rput[t](12.5,33.75){$x''_\ast$}
\rput[t](25,33.75){$x''_\bullet$}
\end{pspicture}

\end{minipage}
\begin{minipage}[c]{.79\textwidth}
 Applying the map $d_{\beta-1}^\epsilon$, we get the result pictured on the left.
 The new configuration represents the action of $d^{\epsilon}_{\beta-1} \circ d^{\epsilon}_{\alpha}$ on $(x_1,\dots,x_k)$ and is labeled $(x''_1,\dots,x''_k)$ on the picture.
\end{minipage}\\[.5em]
Starting from the configuration $(x_1,\dots,x_k)$, we can also send $x_{i_\beta}$ first in $A$ and after send $x_{i_\alpha}$. This leads to the action of the map $d^{\epsilon}_{\alpha}\circ d^{\epsilon}_\beta$ on the configuration $(x_1,\dots,x_k)$. We can observe that 
$d^{\epsilon}_{\beta-1} \circ d^{\epsilon}_{\alpha}$ and $d^{\epsilon}_{\alpha}\circ d^{\epsilon}_\beta$ do not agree exactly but up to homotopy. The homotopy consists of sliding the points $x''_{i_\alpha}$ and  $x''_{i_\beta}$ along the retraction $r$ as sketched below with dotted arrows.
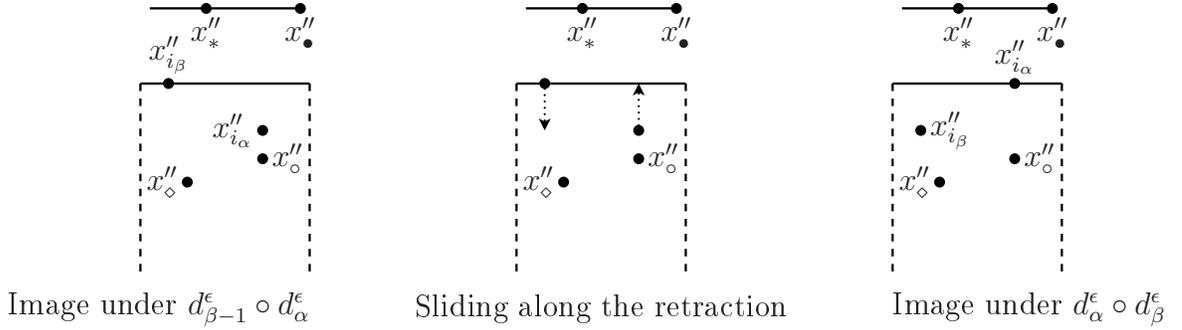
\begin{center}
\psset{xunit=1mm,yunit=1mm,runit=1mm}
\psset{linewidth=0.3,dotsep=1,hatchwidth=0.3,hatchsep=1.5,shadowsize=1}
\psset{dotsize=0.7 2.5,dotscale=1 1,fillcolor=black}
\psset{arrowsize=1 2,arrowlength=1,arrowinset=0.25,tbarsize=0.7 5,bracketlength=0.15,rbracketlength=0.15}
\begin{pspicture}(0,0)(126.25,45)
\psline(3.75,35)(26.25,35)
\psline(5,45)(25,45)
\psdots[](7.5,35)
(20,28.75)
\psdots[](12.5,45)
(25,45)
\psdots[](10,21.88)
(20,25)
\psline[linestyle=dashed,dash=1 1](26.25,35)(26.25,10)
\psline[linestyle=dashed,dash=1 1](3.75,35)(3.75,10)
\rput[l](21.25,25){$x''_\circ$}
\rput[r](8.75,21.88){$x''_\diamond$}
\rput[b](7.5,36.25){$x''_{i_\beta}$}
\rput[r](18.75,28.75){$x''_{i_\alpha}$}
\rput[t](12.5,43.75){$x''_\ast$}
\rput[t](25,43.75){$x''_\bullet$}
\psline(53.75,35)(76.25,35)
\psline(55,45)(75,45)
\psdots[](57.5,35)
(70,28.75)
\psdots[](62.5,45)
(75,45)
\psdots[](60,21.88)
(70,25)
\psline[linestyle=dashed,dash=1 1](76.25,35)(76.25,10)
\psline[linestyle=dashed,dash=1 1](53.75,35)(53.75,10)
\rput[l](71.25,25){$x''_\circ$}
\rput[r](58.75,21.88){$x''_\diamond$}
\rput[t](62.5,43.75){$x''_\ast$}
\rput[t](75,43.75){$x''_\bullet$}
\psline(103.75,35)(126.25,35)
\psline(105,45)(125,45)
\psdots[](107.5,28.75)
(120,35)
\psdots[](112.5,45)
(125,45)
\psdots[](110,21.88)
(120,25)
\psline[linestyle=dashed,dash=1 1](126.25,35)(126.25,10)
\psline[linestyle=dashed,dash=1 1](103.75,35)(103.75,10)
\rput[l](121.25,25){$x''_\circ$}
\rput[r](108.75,21.88){$x''_\diamond$}
\rput[l](108.75,28.75){$x''_{i_\beta}$}
\rput[b](120,36.25){$x''_{i_\alpha}$}
\rput[t](112.5,43.75){$x''_\ast$}
\rput[t](125,43.75){$x''_\bullet$}
\psline[linestyle=dotted,dotsep=0.6]{->}(57.5,35)(57.5,28.75)
\psline[linestyle=dotted,dotsep=0.6]{->}(70,28.75)(70,35)
\rput[r](26.25,5){Image under $d^{\epsilon}_{\beta-1} \circ d^{\epsilon}_{\alpha}$}
\rput[l](103.75,5){Image under $d^{\epsilon}_{\alpha} \circ d^{\epsilon}_{\beta}$}
\rput(65,5){Sliding along the retraction}
\end{pspicture}

\end{center}
Observe that for both compositions, the point $x''_{i_\alpha}$ is of the form
$(i^\epsilon(x_{i_\alpha}),t_\alpha)\in\partial A\times [0,2[\subset A$
and similarly $x''_{i_\beta}$ is of the form $(i^\epsilon(x_{i_\beta}),t_\beta)$.
The preceding homotopy consists in modifying the $t_*$ components of those points.
It is well defined since it is assumed that $i^\epsilon(x_{i_\alpha})\neq i^\epsilon(x_{i_\beta})$.
As a consequence, we have a homotopy between the two configurations.
\begin{prop} \label{relation-cubique}
If $\alpha <\beta$ and $(\epsilon',\epsilon)\in\{-1,+1\}^2$, then
$d^{\epsilon}_{\alpha}\circ d^{\epsilon'}_\beta \simeq d^{\epsilon'}_{\beta-1} \circ d^{\epsilon}_{\alpha}$.
\end{prop}
A more general statement is proved in next section.
It provides the coherences between the preceding homotopies
required to give a structure of $h\cube$-space to $X^k_\bullet$.

\subsection{$h\cube$-space associated to an embedding torus}\label{CubeResol}
In this paragraph, we detail the $h\cube$-space structure of $X^k_\bullet$.
More precisely, we exhibit maps satisfying the relations appearing in Definition \ref{def:hcdiagram}.
\begin{PROP} \label{struct-cubique}
$X_\bullet$ admits a $h\cube$-space structure.
\end{PROP}
\begin{DEMO}
We now prove that $X_\bullet$ admits a $h\cube$-space structure 
as defined in Section \ref{sec:CubicalSpaces}.
To prove this result, we introduce maps $b_*^*$.
Let $i=(i_1,\dots,i_q)\in\Cal C_{q}^{p}$ and $\epsilon=(\epsilon_1,\dots,\epsilon_q)\in\{-1,1\}^q$, $q<p$.
The map $b_{i}^\epsilon$ is defined as: 
\[
b_{  i}^\epsilon:
\begin{array}[t]{rcl}
{\Cal C}^k_p\times\FF_p(C)\times\FF_{k-p}(\In A) &\to& {\Cal C}^k_{p-q}\times\FF_{p-q}(C)\times\FF_{k-p+q}(A)\cr
(  m,(c_1,\dots,c_p),(a_1,\dots,a_{k-p})) &\mapsto& ((\complement  i) \wedge  m,(c_{n_1},\dots,c_{n_{p-q}}),(a'_1,\dots,a'_{k-p+q}))\cr
\end{array}\]
\[\text{with }
  n=(n_1,\dots,n_{p-q})=\complement  i\in\Cal C^p_{p-q},\quad
  l=[i,\complement m]\in \Cal C^{k-p+q}_{k-p} ,\quad
  s=\complement  l\in \Cal C^{k-p+q}_{q},
\] 
\[
\quad a'_{l_\alpha}=a_\alpha \quad\text{ and }\quad a'_{s_\beta}=i^{\epsilon_\beta}(c_{i_\beta}) ,\quad \alpha\in\ito{k-p},\beta\in\ito{q}.
\]
Observe that we cannot compose two maps $b_i^\epsilon$ ; for doing that, we need to use the retraction $r=R(\bullet,1)$ defined before (\ref{retract}).
For instance, the map $d_\alpha^\varepsilon$ introduced in Definition \ref{def-dalpha} 
agrees with the composite $b^{(\varepsilon)}_{(\alpha)}\circ R(\bullet,1)$. In that case, $q=1$.
For a general $q\geq 1$, the map $b_i^\epsilon\circ R(\bullet,1)$ has a geometrical description analogous to the one given in Section \ref{Cube2resol}, but sends $q$ points of $C$ into $A$ instead of just $1$.
We now detail the composite with the general retraction $R(\bullet,1-t)$.

Let $t\in ]0,1]$, $ j\in\Cal C_r^{p-q}$ and $\epsilon'=(\epsilon'_1,\dots,\epsilon'_r)\in\{-1,1\}^r$.
Consider the composition of maps $b_{  j}^{\epsilon'} \circ R(\bullet,1-t) \circ b_{  i}^{\epsilon}$, a straightforward calculus shows that this composition is described as follows:
\[
\begin{array}[t]{rcl}
{\Cal C}^k_p\times\FF_p(C)\times\FF_{k-p}(\In A) &\to& {\Cal C}^k_{p-q-r}\times\FF_{p-q-r}(C)\times\FF_{k-p+q+r}(A)\cr
(  m,(c_1,\dots,c_p),(a_1,\dots,a_{k-p-r})) &\mapsto& ((\complement  j) \wedge(\complement  i) \wedge  m,(c_{{n''_1}},\dots,c_{{n''_{p-q-r}}}),(a''_1,\dots,a''_{k-p+q+r}))\cr
\end{array}\]
where
\begin{itemize}
	\item $  n'=(n'_1,\dots,n'_{p-q-r})=\complement  j\in\Cal C^{p-q}_{p-q-r}$, $  n''=(n_{n'_1},\dots,n_{n'_{p-q-r}})=n'\wedge n =\complement (j\vee i)\in\Cal C^{p}_{p-q-r}$,
	\item $a''_{l'_\alpha}=R(a'_{\alpha},1-t)$, 
	 with $l'=[j,i\vee \complement m]\in\Cal C^{k-p+q+r}_{k-p+q}$ and $\alpha\in\ito{k-p+q}$, which implies
\[a''_{l'_{l_\alpha}}=R(a_\alpha,1-t) \text{ if } \alpha\in\ito{k-p}
\text{ and }
a''_{l'_{s_\alpha}}=R(i^{\epsilon_\alpha}(c_{i_\alpha}),1-t) \text{ if } \alpha\in\ito{q},\]
	\item $a''_{s'_\alpha}=i^{\epsilon'_\alpha}(c_{n_{j_\alpha}})$, with $s'=\complement l'\in\Cal C^{k-p+q+r}_{r}$ and $\alpha\in\ito{r}$.
\end{itemize}
Observe that the expression of the elements of the form $a''_*$ shows that $b_{j}^{\epsilon'} \circ R(\bullet,1-t) \circ b_{  i}^{\epsilon}$ can be extended to $t=1$ by a continuous map.
In that case, we have:
\begin{itemize}
	\item $a''_{l''_\alpha}=a_\alpha$, for $\alpha\in\ito{k-p}$ and $l''=l\wedge l'=[i,\complement m]\wedge [j,i\vee \complement m] =[j\vee i,\complement m]\in \Cal C^{k-p+q+r}_{k-p}$.
	\item $a''_{s''_\alpha}=i^{\epsilon''_\alpha}(c_{w_\alpha})$
	 with $\begin{array}[t]|{l}.
	  s''=s' + s\wedge l'=s' + s\wedge \complement s'=s\vee s'=\complement(l\wedge  l')=\complement(l''),\cr
	  w=j\wedge n + i = j\wedge\complement i + i = j\vee i,\cr
	  \epsilon''=\epsilon_{j\vee i}.\cr
	 	\end{array}$
\end{itemize}
Since $(\complement  j) \wedge(\complement  i) \wedge  m=\complement  (j \vee i) \wedge  m$, we recognize exactly the map $b_{  j\vee i}^{\epsilon_{j\vee i}}$.
Finally, for each element
\[(f_{\underline{i}_n}^{\underline \epsilon_n},t_n,f_{\underline{i}_{n-1}}^{\underline \epsilon_{n-1}},t_{n-1},\dots,t_2,f_{\underline{i}_1}^{\underline \epsilon_1};x)\in (\Mor_\cube)_n(p,p-q)\times I^{n-1}\times ({\Cal C}^k_p\times\FF_p(C)\times\FF_{k-p}(A))\]
define  $\begin{array}[t]{l}X_{p-q}(f_{\underline{i}_n}^{\underline \epsilon_n},t_n,f_{\underline{i}_{n-1}}^{\underline \epsilon_{n-1}},t_{n-1},\dots,t_2,f_{\underline{i}_1}^{\underline \epsilon_1};x)=\cr
\qquad   b_{\underline{i}_{n}}^{\underline \epsilon_{n}}\circ R(\bullet,1-t_{n})\circ b_{\underline{i}_{n-1}}^{\underline \epsilon_{n-1}}\circ R(\bullet,1-t_{n-1})\circ\cdots\circ R(\bullet,1-t_2) \circ b_{\underline{i}_1}^{\underline \epsilon_1}\circ R(x,1) .\cr
\end{array}$\\
The previous calculus of $b_{  j}^{\epsilon'} \circ R(\bullet,1-t) \circ b_{  i}^{\epsilon}$ shows that the maps $X_*$ are well defined, give a coherent system of homotopy, and thus define a structure of $h\cube$-space on $X_\bullet$.
\end{DEMO}

Let $\iota$ denotes the inclusion of $A$ into $M$.
Since $\iota$ is an injective map, it extends to a map, labeled $d_0$, from $\FF_k( A)$ into $\FF_k (M)$.
\begin{PROP}
There exists a homotopy between $d_0 \circ d_1^+$ and $d_0 \circ d_1^-$ from
$\Cal C^{k}_{1}\times C\times\FF_{k-1}(A)$ to $\FF_k(M)$.
Hence, the $h\cube$-space $X^k_\bullet$ is augmented by $d_0$ and there is a map $\chi:|X^k_\bullet|\to\FF_k(M)$ induced by this augmentation.
\end{PROP}
In the next section, we will prove the main Theorem by showing that the map $\chi$ is a homotopy equivalence.
\begin{DEMO}
Let $g$ be the inclusion of $C\times[-1,1]$ into $M$.
If $(i,c,(a_1,\dots,a_{k-1}))\in\Cal C^{k}_{1}\times C \times\FF_{k-1}(A) $, then
\[d_0\circ d_1^\epsilon (i,c,(a_1,\dots,a_{k-1}))=d_0(r(a_1),\dots,r(a_{i-1}), i^\epsilon(c),r(a_{i}),\dots,r(a_{k-1})).\]
Define the map $H : \Cal C^{k}_{1}\times C\times \FF_{k-1}(A)\times[-1,1]\to \FF_k(M)$ that sends $(i,c,(a_1,\dots,a_{k-1}),t)$ to
\[\left(\iota\circ R(a_1,|t|),\dots,\iota\circ R(a_{i-1},|t|), g(c,t),\iota\circ R(a_{i},|t|),\dots,\iota\circ R(a_{k-1},|t|)\right).\]
This is a homotopy between 
$H(\bullet,-1)=d_0\circ d_1^-$ and $H(\bullet,+1)=d_0\circ d_1^+$.
Still, we have to check that $H$ is well defined.
For that, we observe:
\begin{itemize}
    \item The elements $\iota\circ R(a_{*},|t|)$ are all distinct since
$\iota\circ R(\bullet,|t|)$ is injective by construction.
	\item The elements $\iota\circ R(a_{*},|t|)$ are in $A$ and the element $g(c,t)$ is in $C\times[-1,1]$. The intersection $C\times[-1,1]\cap A$ is always included in $\partial A$. But now, we have
\[g(c,t)\in\partial A \Leftrightarrow t=\pm 1
\quad\text{and}\quad \iota(R(a_{i}),|t|)\in \partial A \Leftrightarrow t=0.
\]
therefore, the point $g(c,t)\in C\times[-1,1]$ is always distinct of the points
$\iota\circ R(a_{*},|t|)$. 
\end{itemize}
\end{DEMO}

\subsection{Proof of Main Theorem}\label{DEMOTHA}

The proof relies on the Fadell-Neuwirth fibration \cite{Fadell1962}:
for a manifold $A$ without boundaries, the projection  $\pi_{k,p}$ of the first $p$ components of ${\FF_k}(A)$ onto ${\FF_p}(A)$, $p<k$, is a fiber bundle with fiber over $(a_1,\dots,a_p)\in\FF_p(A)$ the space $\FF_{k-p}(A\setminus\{a_1,\dots,a_p\})$.
Observe that this statement becomes false if $A$ is a manifold with boundary.
The following statement solves this problem by replacing the fiber by an homotopy fiber. It also justifies that the configuration spaces of a manifold and of its interior have the same homotopy type.
\begin{PROP}
Let $A$ be a manifold of interior $\In A$.
Denote by $r : A\to A\setminus(\partial A\times [0,1[)$ the retraction introduced in Proposition \ref{retract}.
The natural inclusions ${\FF_k}(\In  A)\stackrel{i}{\hookrightarrow}{\FF_k}(A)$
are fiber homotopy equivalences. In particular, the homotopy fiber of $\FF_k(A)\to \FF_p(A)$ over an element $(a_1,\dots,a_p)\in\FF_p(A)$ is $\FF_{k-p}(A\setminus\{r(a_1),\dots,r(a_p)\})$.
\end{PROP}
\begin{DEMO}
Observe that the following square is commutative.
$$\xymatrix{
{\FF_{k}}(\In A)\ar@{^(->}[r]\ar[d]	& 	{{\FF_{k}}}(A)\ar[d]	\\
{\FF_{p}}(\In A)\ar@{^(->}[r]	& 	{{\FF_{p}}}(A)	\\}$$
Furthermore, the map $R:A\times [0,2]\to A$ defined in \ref{retract} shows that the composition of two successive maps in the line below is homotopic to the identity.
\[\xymatrix{
{\FF_{k}}(\In A)\ar@{^(->}[rr]^i\ar@/^5ex/[rrrr]^{R(\bullet,1)\he R(\bullet,0)=Id_{{\FF_{k}}(\In A)}} &&
{{\FF_{k}}}(A)\ar[rr]^{R(\bullet,1)}\ar@/_5ex/[rrrr]_{R(\bullet,1)\he R(\bullet,0)=Id_{{{\FF_{k}}}(A)}}		&&	{\FF_{k}}(\In A)\ar@{^(->}[rr]^i	&& 	{{\FF_{k}}}(A)	\\}\]
This proves that ${\FF_{k}}(\In A)\hookrightarrow{\FF_{k}}(A)$ is a fiber homotopy equivalence.
\end{DEMO}
We prove now the Main Theorem and precise the maps between the geometrical realization $|X^k_\bullet|$ and $\FF_k(M)$.
\begin{THEO}\label{TH-resolutioncubique}
The map $\chi$, induced by the augmentation $d_0$, is a homotopy equivalence between the geometrical realization $\left|X^k_\bullet\right|$ and the configuration space $\FF_k(M)$:
\[
\hocolim_{\cube_k} \FF_k\left(C\bigsqcup A\right) \he \FF_k(M)=\FF_k\left(\hocolim (\xymatrix{C\ar@<-.2em>[r]_{i^+}\ar@<.2em>[r]^{i^-}&A}) \right).
\]
\end{THEO}
\begin{DEMO}
In order to prove this theorem, we make an induction on $k$.
If $k=1$, the result is trivial. Now, assume that the geometrical realization of $X^{k-1}_\bullet$ has the same homotopy type as $\FF_{k-1}(M)$.
In order to apply Proposition  \ref{DoubleCubiq}, observe that
$
\Cal C^k_{p}=\left(\Cal C^{k-1}_{p-1}+\{k\}\right) \bigsqcup \Cal C^{k-1}_{p}.
$
Hence, we have the equality:
$$
\Cal C^k_{p}\times{\FF}_{p}(C)\times {\FF}_{k-p}(A)=
\Cal C^{k-1}_{p-1}\times{\FF}_{p}(C)\times {\FF}_{k-p}(A) \bigsqcup \Cal C^{k-1}_{p}\times{\FF}_{p}(C)\times {\FF}_{k-p}(A).
$$
Consider the following diagram:
$$\xymatrix@R=1.1em{
{\Cal C}^{k-1}_{0}\times{\FF}_{1}(C)\times {\FF}_{k-1}(A) \ar@<-.2em>[rrr]_{d_{1}^+}\ar@<.2em>[rrr]^{d_{1}^-}
&&&{\Cal C}^{k-1}_{0}\times{\FF}_{k}(A)\\
{\Cal C}^{k-1}_{1}\times{\FF}_{2}(C)\times {\FF}_{k-2}(A) \ar@<-.2em>[rrr]_{d_{2}^+}\ar@<.2em>[rrr]^{d_{2}^-}\ar@<-.2em>[u]_{d_{1}^-}\ar@<.2em>[u]^{d_{1}^+}
&&&{\Cal C}^{k-1}_{1}\times{\FF}_{1}(C)\times {\FF}_{k-1}(A)
\ar@<-.2em>[u]_{d_{1}^-}\ar@<.2em>[u]^{d_{1}^+}\\
{\Cal C}^{k-1}_{2}\times{\FF}_{3}(C)\times {\FF}_{k-3}(A) \ar@<-.2em>[rrr]_{d_{3}^+}\ar@<.2em>[rrr]^{d_{3}^-}
\ar@<-.5em>[u]_{d_{1}^-}\ar@<.5em>[u]^{d_{2}^+}\ar@3{{}{.}{}}[u]
&&&{\Cal C}^{k-1}_{2}\times{\FF}_{2}(C)\times {\FF}_{k-2}(A)
\ar@<-.5em>[u]_{d_{1}^-}\ar@<.5em>[u]^{d_{2}^+}\ar@3{{}{.}{}}[u]\\
{\vdots}
\ar@{.>}@<-.2em>[rrr]\ar@{.>}@<.2em>[rrr]
\ar@<-.5em>[u]_{d_{1}^-}\ar@<.5em>[u]^{d_{3}^+}\ar@3{{}{.}{}}[u]
&&&{\vdots}\ar@<-.5em>[u]_{d_{1}^-}\ar@<.5em>[u]^{d_{3}^+}\ar@3{{}{.}{}}[u]\\
{\Cal C}^{k-1}_{p-2}\times{\FF}_{p-1}(C)\times {\FF}_{k-p+1}(A) \ar@<-.2em>[rrr]_{d_{p-1}^+}\ar@<.2em>[rrr]^{d_{p-1}^-}
\ar@<-.5em>[u]_{d_{1}^-}\ar@<.5em>[u]^{d_{p-2}^+}\ar@3{{}{.}{}}[u]
&&&{\Cal C}^{k-1}_{p-2}\times{\FF}_{p-2}(C)\times {\FF}_{k-p+2}(A) \ar@<-.5em>[u]_{d_{1}^-}\ar@<.5em>[u]^{d_{p-2}^+}\ar@3{{}{.}{}}[u]\\
{\Cal C}^{k-1}_{p-1}\times{\FF}_{p}(C)\times {\FF}_{k-p}(A) \ar@<-.2em>[rrr]_{d_{p}^+}\ar@<.2em>[rrr]^{d_{p}^-}\ar@<-.5em>[u]_{d_{1}^-}\ar@<.5em>[u]^{d_{p-1}^+}\ar@3{{}{.}{}}[u]
&&&{\Cal C}^{k-1}_{p-1}\times{\FF}_{p-1}(C)\times {\FF}_{k-p+1}(A) \ar@<-.5em>[u]_{d_{1}^-}\ar@<.5em>[u]^{d_{p-1}^+}\ar@3{{}{.}{}}[u]\\
{\vdots}
\ar@{.>}@<-.2em>[rrr]\ar@{.>}@<.2em>[rrr]
\ar@<-.5em>[u]_{d_{1}^-}\ar@<.5em>[u]^{d_{p}^+}\ar@3{{}{.}{}}[u]
&&
&	 {\vdots}\ar@<-.5em>[u]_{d_{1}^-}\ar@<.5em>[u]^{d_{p}^+}\ar@3{{}{.}{}}[u]\\
{\Cal C}^{k-1}_{k-2}\times{\FF}_{k-1}(C)\times {\FF}_{1}(A)	
\ar@<-.2em>[rrr]_{d_{k-1}^+}\ar@<.2em>[rrr]^{d_{k-1}^-}\ar@<-.5em>[u]_{d_{1}^-}\ar@<.5em>[u]^{d_{k-2}^+}\ar@3{{}{.}{}}[u]
&&&	 {\Cal C}^{k-1}_{k-2}\times{\FF}_{k-2}(C)\times {\FF}_{2}(A)\ar@<-.5em>[u]_{d_{1}^-}\ar@<.5em>[u]^{d_{k-2}^+}\ar@3{{}{.}{}}[u]\\
{\FF}_k(C)	
\ar@<-.2em>[rrr]_{d_{k}^+}\ar@<.2em>[rrr]^{d_{k}^-}\ar@<-.5em>[u]_{d_{1}^-}\ar@<.5em>[u]^{d_{k-1}^+}\ar@3{{}{.}{}}[u]
&&&	{\Cal C}^{k-1}_{k-1}\times{\FF}_{k-1}(C)\times {\FF}_1(A) \ar@<-.5em>[u]_{d_{1}^-}\ar@<.5em>[u]^{d_{k-1}^+}\ar@3{{}{.}{}}[u]\\
}$$
Label $C_\bullet^k$ the $h\cube$-space on the left hand side of the diagram and
$A_\bullet^k$ the $h\cube$-space on the right hand side.
Proposition \ref{DoubleCubiq} shows that there exists a $h\cube$-space
$\xymatrix{|C_\bullet^k|\ar@<-.2em>[r]_{d^+}\ar@<.2em>[r]^{d^-}&	|A_\bullet^k|}$
and that its geometrical realization is the same as the one of $X_\bullet^k$.
Furthermore, $|A_\bullet^k|$ is also augmented by the map $d_0:\FF_k(A)\hookrightarrow\FF_k(M)$.
This augmentation induces a map $\chi_A:|A_\bullet^k|\to \FF_k(M)$ that fits in the following homotopy commutative diagram $(\ddag)$:
\[\xymatrix{
\FF_k(A)\ar[r]^-{d_0}	\ar@{^(->}[d]		&			\FF_k(M)		\\
|A_\bullet^k|\ar@{^(->}[r]\ar[ur]_{\chi_A}	&			|X_\bullet^k|\ar[u]_\chi\\
}\]
In order to study $A_\bullet^k$ and $C_\bullet^k$, define
$\pi_k^M:\FF_k(M) \longrightarrow M$ as the projection on the last coordinate of $\FF_k(M)$,
$\pi_k^A:\Cal C^{k-1}_p\times{\FF}_p(C)\times {\FF}_{k-p}(A) \longrightarrow A$ as the projection on the last coordinate of  $\FF_{k-p}(A)$ and 
$\pi_k^C:\left(\Cal C^{k-1}_{p-1}+\{k\}\right)\times{\FF}_{p}(C)\times {\FF}_{k-p}(A) \to C$ as the projection on the last coordinate of  $\FF_{p}(C)$. 
Observe that the $h\cube$-space $A^k_\bullet$ is augmented by the map $\pi_k^A:{\FF}_{k}(A)\to A$. Hence there is an induced map, still denoted $\pi_k^A$, from $|A_\bullet^k|$ to $A$.
Also remark that $\iota\circ \pi_k^A = \pi_k^M \circ d_0$, and consequently the square
$$
\xymatrix{
|A_\bullet^k|\ar[r]^-{\chi_A}\ar[d]^{\pi_k^A}		&	{\FF_k(M)}\ar[d]^{\pi_k^M}\\
A\ar[r]_-{\iota}				&	M\\
}$$
is homotopy commutative.
We claim that this square is a homotopy pullback.
Indeed the map $\chi_A$ restricts to a homotopy equivalence between the homotopy fibers of the vertical maps.
In order to show it, we know from  V.~Puppe \cite{Puppe1974} that the homotopy fiber of $\pi_k^A : |A_\bullet^k|\to A$ is the geometrical realization of the $h\cube$-space obtained by restriction to the homotopy fibers.
More precisely, the homotopy fiber of $\pi_k^A$ over  $a\in A$ is the geometrical realization of the $h\cube$-space below:
$$\xymatrix@R=1.4em{
{\Cal C}^{k-1}_{0}\times {\FF}_{k-1}(A\setminus r(a))\\
{\vdots}\ar@<-.2em>[u]_{d_{1}^-}\ar@<.2em>[u]^{d_{1}^+}\\
{\Cal C}^{k-1}_{p-1}\times{\FF}_{p-1}(C)\times {\FF}_{k-p}(A\setminus r(a))
\ar@<-.5em>[u]_{d_{1}^-}\ar@<.5em>[u]^{d_{p-1}^+}\ar@3{{}{.}{}}[u]\\
{\Cal C}^{k-1}_{p}\times{\FF}_{p}(C)\times {\FF}_{k-p-1}(A\setminus r(a))
\ar@<-.5em>[u]_{d_{1}^-}\ar@<.5em>[u]^{d_{p}^+}\ar@3{{}{.}{}}[u]\\
{\vdots}\ar@<-.5em>[u]_{d_{1}^-}\ar@<.5em>[u]^{d_{p+1}^+}\ar@3{{}{.}{}}[u]\\
{\Cal C}^{k-1}_{k-2}\times{\FF}_{k-2}(C)\times {\FF}_{1}(A\setminus r(a))
\ar@<-.5em>[u]_{d_{1}^-}\ar@<.5em>[u]^{d_{k-2}^+}\ar@3{{}{.}{}}[u]\\
{\Cal C}^{k-1}_{k-1}\times{\FF}_{k-1}(C)
\ar@<-.5em>[u]_{d_{1}^-}\ar@<.5em>[u]^{d_{k-1}^+}\ar@3{{}{.}{}}[u]\\
}$$
By our induction hypothesis, the geometrical realization of this $h\cube$-space has the homotopy type of $\FF_{k-1}(M\setminus r(a))$. Therefore, we get the annouced homotopy pullback.
Using the same argument with $C_\bullet^k$, we show that the outer square of the following diagram is also a homotopy pullback:
$$\xymatrix{
|C_\bullet^k|\ar[r]^{d^\epsilon}\ar[d]^{\pi_k^C}		&|A_\bullet^k|\ar[r]^{\chi_A}\ar[d]^{\pi_k^A}	&{\FF_k(M)}\ar[d]^{\pi_k^M}\\
C\ar[r]_{i^\epsilon}							&				A\ar[r]_\iota				&	M\\
}$$
Furthermore, we have the following equality of maps in $C_\bullet^k$:
\[ d_1^\epsilon\circ\pi_k^C = r\circ i^\epsilon\circ\pi_k^C \he i^\epsilon\circ\pi_k^C = \pi_k^A\circ d_p^\epsilon :
 {\Cal C}^{k-1}_{p-1}\times{\FF}_{p}(C)\times {\FF}_{k-p}(A) \to A.\]
Consequently, the left square is also a homotopy pullback.
Let $F$ be the common homotopy fiber of the three previous vertical arrows.
According to Lemma \ref{DoubleCubiq}, the geometrical realization of the $h\cube$-space $\xymatrix{|C_\bullet^k|\ar@<-0.5ex>[r]_{d^-}\ar@<0.5ex>[r]^{d^+}	& |A_\bullet^k|}$ is the space $|X_\bullet^k|$  and the geometrical realization of $\xymatrix{C\ar@<-0.5ex>[r]_{i^-}\ar@<0.5ex>[r]^{i^+}	& A}$ is $M$.
Using \cite{Puppe1974}, we know that the induced maps between those homotopy colimits $|X^k_\bullet|\to M$ has also $F$ for homotopy fiber. Moreover, since the diagram $(\ddag)$ is homotopy commutative, the square
$$\xymatrix{
|X_\bullet^k|\ar[r]^{\chi}\ar[d]^-{\pi_k}	&{\FF_k(M)}\ar[d]^{\pi_k}\\
 \displaystyle \left|C\mathop{\longrightarrow\atop\longrightarrow}_{i^+}^{i^-} A\right| \ar[r]^-\he &	M\\
}$$
is a homotopy pullback. 
Since the bottom map is a homotopy equivalence, the top map 
$\chi:|X_\bullet^k|\to\FF_k(M)$ is also a homotopy equivalence.
\end{DEMO}

\section{Application to braid groups}\label{sec:Braids}
In this section, we explain the method to compute the pure braid groups of $M$ and detail the case of the Möbius band.

\subsection{Artin braid groups}\label{ArtinBraids}

The braid groups have been introduced by E.~Artin \cite{Artin1947}.
Here we adopt the point of view where braid groups are defined in terms of fundamental group of a configuration space as given in \cite{Fox1962} or, more generally, in \cite{Vershinin1999}.
It is known that $\FF_{k}(\RR^2)$ is an Eilenberg-MacLane space $K(P_k,1)$ where $P_k$ is called the \emph{pure braid group on $k$ strands}.
We extend this definition as follows.
\begin{DEFI}
Let $M$ be a connected manifold of dimension bigger than $2$.
The pure braid group on $k$ strands of $M$ is the group
\[ P_k(M) = \pi_1(\FF_k(M)).\]
\end{DEFI}
In this section, we give a description of the group $P_k$.

Let $\ito{k}$ be the set of integers between $1$ and $k$.
Let $q_1<q_2<\dots<q_k$ be some fixed real numbers.
The image of $q_*$ by the the natural inclusion $\RR\cong\RR\times\{0\}\subset\RR^2$ is still denoted $q_*$.
\emph{The configuration $Q_k=(q_1,\dots,q_k)$ is the base point of the space $\FF_k(\RR^2)$.}
Finally, let $\epsilon=\inf_{i\neq j} |q_j-q_i|/2$ and $\vec{\epsilon}=(\epsilon,0)\in\RR^2$.

For $(i,j)\in\ito{k}^2$, we define a path $t_{j,i}=(\tau_1,\dots,\tau_k):[0,1]\to\FF_k(\RR^2)$
joining the configuration  $(q_1,\dots,q_k)$ to the configuration $(q_1,\dots,q_{j-1},q_i+\vec{\epsilon},q_{j+1},\dots,q_k)$ such that
$\tau_r(t)=q_r$	if $r\neq j$ and $\tau_j(t)\in\RR\times\RR^+$.
Define the loop $\alpha_{j,i}:S^1\to \FF_k(\RR^2)$ by 
$\alpha_{j,i}(\xi)=(q_1,\dots,q_{j-1},q_i+\epsilon\xi,q_{j+1},\dots,q_k)$.
The class $A_{j,i}\in\pi_1(\FF_k(\RR^2),Q_k)$ is defined as the homotopy class
of the loop $t_{j,i} *\alpha_{j,i} * t_{j,i}^{-1}$ where $*$ denote the composition of paths.
Observe that the homotopy class of $A_{j,i}$ is independent of the choice of $t_{j,i}$ and $\epsilon$.

Whenever no confusion is  possible, we denote in the same way the loops and their homotopy classes and we omit to write down the composition operation $*$.

The element $A_{j,i}= t_{j,i} \alpha_{j,i}  t_{j,i}^{-1}$ is represented by the following diagram.
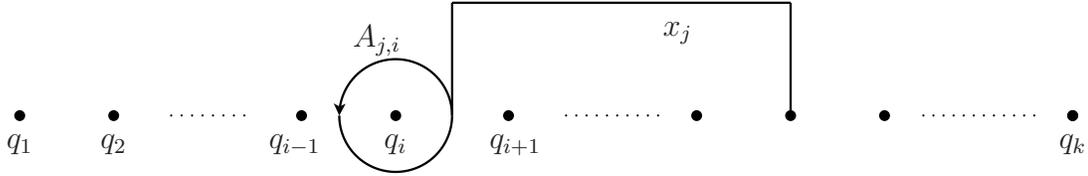
\begin{center}
\psset{xunit=1mm,yunit=1mm,runit=1mm}
\psset{linewidth=0.3,dotsep=1,hatchwidth=0.3,hatchsep=1.5,shadowsize=1}
\psset{dotsize=0.7 2.5,dotscale=1 1,fillcolor=black}
\psset{arrowsize=1 2,arrowlength=1,arrowinset=0.25,tbarsize=0.7 5,bracketlength=0.15,rbracketlength=0.15}
\begin{pspicture}(0,0)(140,22.5)
\psdots[linestyle=none](0,7.5)
(12.5,7.5)
\psdots[linestyle=none](37.5,7.5)
(50,7.5)
\psdots[linestyle=none](90,7.5)
(65,7.5)
\psdots[linestyle=none](102.5,7.5)
(115,7.5)
\psdots[linestyle=none](140,7.5)
(140,7.5)
\rput[t](0,5){$q_1$}
\rput[t](12.5,5){$q_2$}
\rput[t](140,5){$q_k$}
\rput[t](50,5){$q_i$}
\rput[t](87.5,20){$x_j$}
\rput(47.5,17.5){$A_{j,i}$}
\psline[linestyle=dotted](120,7.5)(135,7.5)
\psline[linestyle=dotted](72.5,7.5)(85,7.5)
\psline[linestyle=dotted](20,7.5)(30,7.5)
\rput[tr](40,5){$q_{i-1}$}
\rput[tl](62.5,5){$q_{i+1}$}
\psline(102.5,7.5)(102.5,22.5)
\psline(102.5,22.5)(57.5,22.5)
\psline(57.5,22.5)(57.5,7.5)
\rput{0}(50,7.5){\parametricplot[arrows=->]{0}{180}{ t cos 7.5 mul t sin 7.5 mul }}
\rput{0}(50,7.5){\parametricplot[arrows=-]{180}{360}{ t cos 7.5 mul t sin 7.5 mul }}
\end{pspicture}

\end{center}
The $r$-th particle of the configuration is labeled $x_r$, except when it is fixed at the base point and is denoted $q_r$ in that case.

Analogously, we define the class  $B_{j,i}$ by requesting $\tau_j(t)\in\RR\times\RR^-$. It is represented by the next diagram.
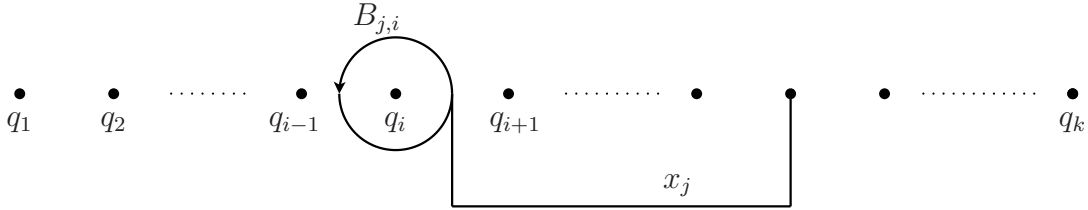
\begin{center}
\psset{xunit=1mm,yunit=1mm,runit=1mm}
\psset{linewidth=0.3,dotsep=1,hatchwidth=0.3,hatchsep=1.5,shadowsize=1}
\psset{dotsize=0.7 2.5,dotscale=1 1,fillcolor=black}
\psset{arrowsize=1 2,arrowlength=1,arrowinset=0.25,tbarsize=0.7 5,bracketlength=0.15,rbracketlength=0.15}
\begin{pspicture}(0,0)(140,25)
\psdots[linestyle=none](0,15)
(12.5,15)
\psdots[linestyle=none](37.5,15)
(50,15)
\psdots[linestyle=none](90,15)
(65,15)
\psdots[linestyle=none](102.5,15)
(115,15)
\psdots[linestyle=none](140,15)
(140,15)
\rput[t](0,12.5){$q_1$}
\rput[t](12.5,12.5){$q_2$}
\rput[t](140,12.5){$q_k$}
\rput[t](50,12.5){$q_i$}
\rput[B](87.5,2.5){$x_j$}
\rput(47.5,25){$B_{j,i}$}
\psline[linestyle=dotted](120,15)(135,15)
\psline[linestyle=dotted](72.5,15)(85,15)
\psline[linestyle=dotted](20,15)(30,15)
\rput[tr](40,12.5){$q_{i-1}$}
\rput[tl](62.5,12.5){$q_{i+1}$}
\psline(102.5,0)(102.5,15)
\psline(102.5,0)(57.5,0)
\psline(57.5,15)(57.5,0)
\rput{0}(50,15){\parametricplot[arrows=->]{0}{180}{ t cos 7.5 mul t sin 7.5 mul }}
\rput{0}(50,15){\parametricplot[arrows=-]{180}{360}{ t cos 7.5 mul t sin 7.5 mul }}
\end{pspicture}

\end{center}
The class $B_{j,i}$ is related to the classes $A_{*,*}$ by the relation 
$$B_{j,i}=A_{j,j-1}^{-1}A_{j,j-2}^{-1}\cdots A_{j,i+1}^{-1} A_{j,i} A_{j,i+1}\cdots A_{j,j-2} A_{j,j-1}.$$

The following theorem gives a description of the group $P_k$.
\begin{THEO}\label{equ-tresses}
The group  $P_k$ admits the following presentation:\\
Generators : $A_{j,i}$ with $1\leq i < j \leq k$.\\
Relations : 
\[\begin{array}{lrl}
	(1)		&	[A_{j,i},A_{r,i}A_{r,j}]=1	&	\text{ if }\  1\leq i < j < r ,\cr
	(2)		&	[A_{r,i},A_{r,j}A_{j,i}]=1	&	\text{ if }\  1\leq i < j < r  ,\cr
	(3)		&	[A_{s,r},A_{j,i}]=1			&	\text{ if }\  1\leq i < j < r < s ,\cr
	(4)		&	[A_{s,i},A_{r,j}]=1			&	\text{ if }\  1\leq i < j < r < s ,\cr
	(5)		&	[A_{s,j},A_{r,j}^{-1}A_{r,i}A_{r,j}]=1	&	\text{ if }\  1\leq i < j < r < s ,\cr
	(6)		&	[A_{s,j},A_{s,r}A_{r,i}A_{s,r}^{-1}]=1	&	\text{ if }\  1\leq i < j < r < s ,\cr
\end{array}\]
where the commutator $[A,B]$ denotes the element $A^{-1}B^{-1}AB$.
\end{THEO}
A complete proof of this theorem can be found in \cite{Fox1962,Fadell2001}.
It is also a consequence of the next lemma that we will use afterward.
\begin{LEMM}\label{commut-lacets}
Let $M$ be a manifold and $Q_k=(q_1,\dots,q_k)\in\FF_k(M)$ be a configuration in $M$.\\
Let $U$ and $V$ be two disjoint subsets of $\ito{k}$.\\
Let $\tau=(\tau_1,\dots,\tau_k)$ be a path in $\FF_k(M)$ starting at $Q_k$ and such that $\tau_u$ is constant if $u\notin U$.\\
Let $\gamma=(\gamma_1,\dots,\gamma_k)$ be a loop in $\FF_k(M)$ pointed at $Q_k$ such that $\gamma_v$ is constant if $v\notin V$.\\
Suppose that for all couple $(u,v)\in U\times V$ and all $(t_1,t_2)\in[0,1]^2$, we have $\tau_u(t_1)\neq\gamma_v(t_2)$.
Then the loop 
$\omega=(\omega_1,\dots,\omega_k):[0,1]\to\FF_k(M)$ defined by 
$\left\{\begin{array}[c]{ll}
\omega_u=\tau_u * 1 * \tau_u^{-1} & \text{if } u\in U\cr
\omega_v=1 *	\gamma_v * 1 & \text{if } v\in V	\cr
\end{array}\right.$ is well defined and is homotopic to the path $\gamma$.
In particular, if $\gamma$ is a loop, we have $\omega=\tau\gamma\tau^{-1}$ and so $[\gamma,\tau]=1 \in \pi_1(\FF_k(M))$.
\end{LEMM}
\begin{DEMO}
Let $s\in[0,1]$ and $\tau_{u,s}$ be the path in $M$ defined by $\tau_{u,s}(t)=\tau_{u}(st)$.
Define a loop $\Omega_s=(\omega_{1,s},\dots,\omega_{k,s}):[0,1]\to\FF_k(M)$  by
$\left\{\begin{array}[c]{ll}
\omega_{u,s}=\tau_{u,s} * 1 * \tau_{u,s}^{-1} & \text{if } u\in U\cr
\omega_{v,s}=1 *	\gamma_v * 1 & \text{if } v\in V	\cr
\end{array}\right.$ where the path named  $1$ is the constant path at the right point. The path $\Omega_s$ is well defined because for every couple $(u,v)\in U\times V$ and all $t\in[0,1]$, we have $\tau_{u,s}(t)=\tau_u(st)\neq\gamma_v(t)$.
Therefore, the loops $\omega=\Omega_1$ and $\Omega_0=1*\gamma*1\he\gamma$ are in the same homotopy class.
\end{DEMO}
In the case of $3$ or $4$ particles, the proof of Theorem \ref{equ-tresses} is contained in the following diagrams.
\begin{center}
\psset{xunit=1mm,yunit=1mm,runit=1mm}
\psset{linewidth=0.3,dotsep=1,hatchwidth=0.3,hatchsep=1.5,shadowsize=1}
\psset{dotsize=0.7 2.5,dotscale=1 1,fillcolor=black}
\psset{arrowsize=1 2,arrowlength=1,arrowinset=0.25,tbarsize=0.7 5,bracketlength=0.15,rbracketlength=0.15}
\begin{pspicture}(0,0)(155,81.25)
\psdots[](90,40)
(110,40)
\psdots[](130,40)
(150,40)
\pscustom[]{\psbezier(150,40)(150,50)(140,50)(130,50)
\psbezier(130,50)(120,50)(95,50)(90,50)
\psbezier(90,50)(85,50)(80,45)(80,40)
\psbezier(80,40)(80,35)(85,30)(90,30)
\psbezier(90,30)(95,30)(100,35)(100,40)
\psbezier(100,40)(100,45)(95,50)(90,50)
}
\pscustom[]{\psbezier(130,40)(130,45)(125,45)(120,45)
\psbezier(120,45)(115,45)(115,45)(110,45)
\psbezier(110,45)(105,45)(105,40)(105,40)
\psbezier(105,40)(105,40)(105,35)(110,35)
\psbezier(110,35)(115,35)(115,40)(115,40)
\psbezier(115,40)(115,40)(115,45)(110,45)
}
\psdots[](10,40)
(30,40)
\psdots[](50,40)
(70,40)
\pscustom[]{\psbezier(30,40)(30,45)(30,50)(20,50)
\psbezier(20,50)(10,50)(15,50)(10,50)
\psbezier(10,50)(5,50)(0,45)(0,40)
\psbezier(0,40)(0,35)(5,30)(10,30)
\psbezier(10,30)(15,30)(20,35)(20,40)
\psbezier(20,40)(20,45)(15,50)(10,50)
}
\pscustom[]{\psbezier(70,40)(70,45)(70,50)(60,50)
\psbezier(60,50)(50,50)(55,50)(50,50)
\psbezier(50,50)(45,50)(40,45)(40,40)
\psbezier(40,40)(40,35)(45,30)(50,30)
\psbezier(50,30)(55,30)(60,35)(60,40)
\psbezier(60,40)(60,45)(55,50)(50,50)
}
\psdots[](90,10)
(110,10)
\psdots[](130,10)
(150,10)
\pscustom[]{\psbezier(130,10)(130,15)(130,20)(112.5,20)
\psbezier(112.5,20)(95,20)(95,15)(90,15)
\psbezier(90,15)(85,15)(85,10)(85,10)
\psbezier(85,10)(85,10)(85,5)(90,5)
\psbezier(90,5)(95,5)(95,10)(95,10)
\psbezier(95,10)(95,10)(95,15)(90,15)
}
\pscustom[]{\psbezier(150,10)(150,15)(150,20)(145,20)
\psbezier(145,20)(140,20)(135,20)(135,10)
\psbezier(135,10)(135,0)(125,0)(125,10)
\psbezier(125,10)(125,20)(113.75,17.5)(110,17.5)
\psbezier(110,17.5)(106.25,17.5)(102.5,13.75)(102.5,10)
\psbezier(102.5,10)(102.5,6.25)(105,2.5)(110,2.5)
\psbezier(110,2.5)(115,2.5)(117.5,6.25)(117.5,10)
\psbezier(117.5,10)(117.5,13.75)(113.75,17.5)(110,17.5)
}
\psdots[](10,10)
(30,10)
\psdots[](50,10)
(70,10)
\pscustom[]{\psbezier(70,10)(70,15)(70,20)(60,20)
\psbezier(60,20)(50,20)(35,15)(30,15)
\psbezier(30,15)(25,15)(25,10)(25,10)
\psbezier(25,10)(25,10)(25,5)(30,5)
\psbezier(30,5)(35,5)(35,10)(35,10)
\psbezier(35,10)(35,10)(35,15)(30,15)
}
\pscustom[]{\psbezier(50,10)(50,15)(50,15)(45,15)
\psbezier(45,15)(40,15)(40,15)(37.5,10)
\psbezier(37.5,10)(35,5)(35,2.5)(30,2.5)
\psbezier(30,2.5)(25,2.5)(25,5)(22.5,10)
\psbezier(22.5,10)(20,15)(15,15)(10,15)
\psbezier(10,15)(5,15)(5,10)(5,10)
\psbezier(5,10)(5,10)(5,5)(10,5)
\psbezier(10,5)(15,5)(15,10)(15,10)
\psbezier(15,10)(15,10)(15,15)(10,15)
}
\psdots[](10,70)
(10,70)
\pscustom[]{\psbezier(40,70)(40,80)(30,77.5)(22.5,77.5)
\psbezier(22.5,77.5)(15,77.5)(15,77.5)(10,77.5)
\psbezier(10,77.5)(5,77.5)(2.5,75)(2.5,70)
\psbezier(2.5,70)(2.5,65)(5,62.5)(10,62.5)
\psbezier(10,62.5)(15,62.5)(17.5,65)(17.5,70)
\psbezier(17.5,70)(17.5,75)(15,77.5)(10,77.5)
}
\pscustom[]{\psbezier(70,70)(70,80)(63.75,80)(56.25,80)
\psbezier(56.25,80)(48.75,80)(20,80)(10,80)
\psbezier(10,80)(0,80)(0,75)(0,70)
\psbezier(0,70)(0,65)(2.5,60)(10,60)
\psbezier(10,60)(17.5,60)(20,65)(20,70)
\psbezier(20,70)(20,75)(30,75)(35,75)
\psbezier(35,75)(40,75)(30,61.25)(40,61.25)
\psbezier(40,61.25)(50,61.25)(40,75)(45,75)
\psbezier(45,75)(50,75)(56.25,75)(60,75)
\psbezier(60,75)(63.75,75)(70,75)(70,70)
}
\psbezier(135,70)(135,70)(135,70)(135,70)
\psdots[](70,70)
(40,70)
\psdots[](90,70)
(90,70)
\pscustom[]{\psbezier(150,70)(150,77.5)(150,77.5)(142.5,77.5)
\psbezier(142.5,77.5)(135,77.5)(95,77.5)(90,77.5)
\psbezier(90,77.5)(85,77.5)(82.5,75)(82.5,70)
\psbezier(82.5,70)(82.5,65)(85,62.5)(90,62.5)
\psbezier(90,62.5)(95,62.5)(97.5,65)(97.5,70)
\psbezier(97.5,70)(97.5,75)(95,77.5)(90,77.5)
}
\psdots[](150,70)
(120,70)
\pscustom[]{\psbezier(120,70)(120,77.5)(126.25,76.25)(135,76.25)
\psbezier(135,76.25)(143.75,76.25)(146.25,76.25)(146.25,70)
\psbezier(146.25,70)(146.25,63.75)(146.25,60)(150,60)
\psbezier(150,60)(153.75,60)(155,63.75)(155,70)
\psbezier(155,70)(155,76.25)(155,81.25)(150,81.25)
\psbezier(150,81.25)(145,81.25)(96.25,81.25)(90,81.25)
\psbezier(90,81.25)(83.75,81.25)(80,76.25)(80,70)
\psbezier(80,70)(80,63.75)(82.5,60)(90,60)
\psbezier(90,60)(97.5,60)(100,63.75)(100,70)
\psbezier(100,70)(100,76.25)(103.75,76.25)(110,76.25)
\psbezier(110,76.25)(116.25,76.25)(120,76.25)(120,70)
}
\rput[t](40,60){$[A_{2,1},A_{3,1}A_{3,2}]=1$}
\rput[t](120,60){$[A_{3,1},A_{3,2}A_{2,1}]=1$}
\rput[t](40,30){$[A_{4,3},A_{2,1}]=1$}
\rput[t](120,31.25){$[A_{4,1},A_{3,2}]=1$}
\rput[t](40,0){$[A_{4,2},A_{3,2}^{-1}A_{3,1}A_{3,2}]=1$}
\rput[t](118.75,0){$[A_{3,1},A_{4,3}A_{4,2}A_{4,3}^{-1}]=1$}
\rput[t](10,68.75){$q_1$}
\rput[t](70,68.75){$x_3$}
\rput[t](40,68.75){$x_2$}
\rput[t](90,68.75){$q_1$}
\rput[t](150,68.75){$q_3$}
\rput[t](120,68.75){$q_2$}
\rput[t](10,38.75){$q_1$}
\rput[t](50,38.75){$q_3$}
\rput[t](30,38.75){$x_2$}
\rput[t](70,38.75){$x_4$}
\rput[t](10,8.75){$q_1$}
\rput[t](50,8.75){$x_3$}
\rput[t](30,8.75){$q_2$}
\rput[t](70,8.75){$x_4$}
\rput[t](90,38.75){$q_1$}
\rput[t](130,38.75){$x_3$}
\rput[t](110,38.75){$q_2$}
\rput[t](150,38.75){$x_4$}
\rput[t](90,8.75){$q_1$}
\rput[t](130,8.75){$x_3$}
\rput[t](110,8.75){$q_2$}
\rput[t](150,8.75){$x_4$}
\end{pspicture}
\\*[1.2em]
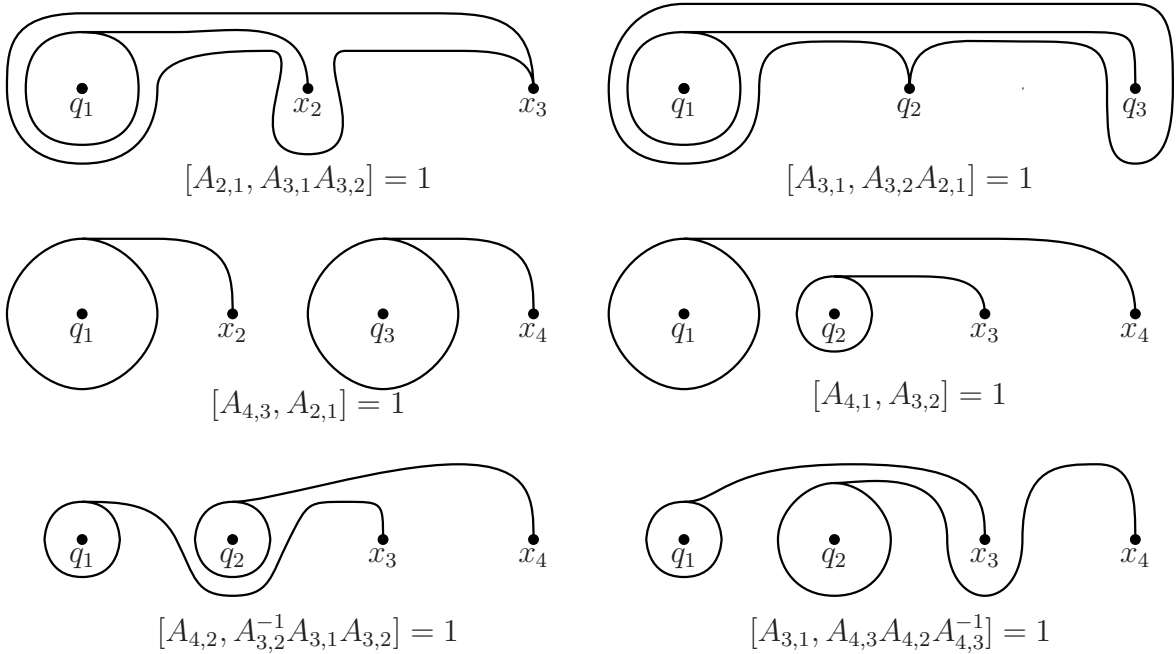
	Figure: Yang-Baxter relations for $3$ or $4$ particles.
\end{center}

\subsection{Practical determination of $P_k(M)$}\label{Method}

As pointed out in Section \ref{sec:CubicalSpaces}, one does not need to know the rigorous structure of a $h\cube$-space in order to calculate its fundamental group.
To calculate the fundamental group of $|X^k_\bullet|$, which is $P_k(M)$, the truncated realization $|X^k_{\leq2}|$ is enough.
In fact, if $(X^k_2)^{(0)}$ denotes a $0$-skeleton for $X^k_2$, Proposition \ref{PI1-cubiq} asserts that the geometrical realization of the $h\cube$-space
\[
Y_\bullet : 
\xymatrix{
(X^k_2)^{(0)} \ar@<-0.5ex>[r]\ar@<0.5ex>[r]\ar@<-1.5ex>[r]_-{d_2^+}\ar@<1.5ex>[r]^-{d_1^-}	&
X^k_1\ar@<-0.5ex>[r]_-{d_1^+}\ar@<0.5ex>[r]^-{d_1^-}&
X_0^k\\
}\]
has the same fundamental group as $|X^k_\bullet|$ which is $P_k(M)$ by the Main Theorem.
Here, we detail how this group can be described by generators and relations.
 Let $Y_\bullet$ be as above and denote by $Q_k\in\FF_k(A)=Y_0$ the common base point of $|Y_{\leq1}|$ and $|Y_{\leq2}|$.

\noindent\textbf{First step: } The space $|Y_{\leq1}|$ is the mapping torus
\[|Y_{\leq1}|=
\hocolim 
\left(\xymatrix{
\ito{k}\times C \times \FF_{k-1}(A)\ar@<-0.5ex>[r]_-{d_1^+}\ar@<0.5ex>[r]^-{d_1^-}&
\FF_k(A)}\right).\]
In order to find the fundamental group of $|Y_{\leq1}|$, 
we apply a Van Kampen like theorem inductively for each path-component of $Y_1=\ito{k}\times C \times \FF_{k-1}(A)$.
If a path component of $Y_1=\ito{k}\times C \times \FF_{k-1}(A)$ is sent to two different path-components of $Y_0=\FF_k(A)$ by the maps $d_1^-$ and $d_1^+$, we can use the usual Van Kampen theorem (see \cite[Theorem 1.20]{Hatcher2002}). In the other case, we use the following variation which can be easily proved.
\begin{prop}[Van Kampen Theorem for a mapping torus]\label{VK-tore}
Let $f^-$ and $f^+$ be two maps between two arcwise connected spaces $Z$ and $Y$.
Let $X$ be the homotopy colimit of the diagram $$\xymatrix{Z\ar@<-.2em>[r]_{f^+}\ar@<.2em>[r]^{f^-}&Y}$$
i.e.
$$X={Y\bigsqcup Z\times[-1,1] \over f^+(z)\sim(z,+1),f^-(z)\sim(z,-1)}$$
Fix two points $y_0\in Y$ and $z_0\in Z$,
a path $y^+$ in $Y$ from $y_0$ to $f^+(z_0)$,
a path $y^-$ in $Y$ from $y_0$ to $f^-(z_0)$,
and $z$ the path in $X$ with support $\{z_0\}\times [-1,1]$.
Then, there exists an isomorphism
$$
\pi_1(X,y_0)\approx {\pi_1(Y,y_0)\star <\rho> \over
 {y^+}f^+(\omega)(y^+)^{-1}\sim \rho^{-1}{y^-}f^-(\omega)(y^-)^{-1}\rho \text{ for } \omega\in\pi_1(Z,z_0)}
$$
where $\star$ denotes the free product of groups and $\rho\in\pi_1(X,y_0)$ is the loop $y^- z (y^+)^{-1}$.
\end{prop}

\noindent\textbf{Second step: }
In order to obtain $|Y_{\leq 2}|$, we glue a $2$-cell on top of $|Y_{\leq1}|$ for each point $x\in Y_2={\left({\Cal C}_2^k\times{\FF_2(C)}\times \FF_{k-2}(A)\right)^{(0)}}$.
This $2$-cell is attached along the loop $\Phi_2(x,\bullet): S^1\to |Y_{\leq1}|$ described in \ref{realization-cube2}.
Also, for each point $x\in Y_2={\left({\Cal C}_2^k\times{\FF_2(C)}\times \FF_{k-2}(A)\right)^{(0)}}$, choose $\alpha_x$ a path in $|X_{\leq 1}|$ from $\Phi_2(x,1)$ to $Q_k$.
Finally, as it is well known (see \cite[Proposition 1.26]{Hatcher2002}), the group $P_k(M)=\pi_1(|X^k_\bullet|,Q_k)=\pi_1(|Y_\bullet|,Q_k)$ admits for presentation 
\[ \pi_1(|X_{\leq 1}|,Q_k)\big/ <\alpha_x  \Phi_2(x,\bullet)  \alpha_x^{-1} \ \big|\  x\in Y_2 >.\]

\subsection{Braids on the Möbius band $\Cal M$}\label{Mobius}

The Möbius band $\Cal M$ can be represented as an embedding torus with $C$ the interval $[-1,1]$ and $A$ the square $[-1,1]\times[-1,1]$.
The maps $i^\epsilon$ are defined by $i^\epsilon(x)=(\epsilon x,\epsilon)$.

We apply the method described in the previous section to $\Cal M$. Observe that it is sufficient to describe the group $P_k(\Cal M)$  to know the homotopy type of $\FF_k(\Cal M)$ since this space is a $K(P_k(\Cal M),1)$.
Indeed, it follows from the long exact sequence of the fibration
$$\bigvee_1^k S^1\he\Cal M \setminus \{q_1,\dots,q_{k-1}\}\to \FF_k(\Cal M)\to \FF_{k-1}(\Cal M)$$
and a trivial induction.

\begin{theo}\label{Present-Mobius}
The group $P_k(\Cal M)$ admits the following presentation:
\begin{itemize}
\item Generators: $\rho_i$ for $1\leq i \leq k$.
\item Relations:
\[
\begin{array}{l>{\text{if }}l}
[\rho_i^{-1},\rho_j^{-1}]=[\rho_{j-1},\rho_j]\cdots [\rho_{i+1},\rho_j]
[\rho_{j},\rho_i][\rho_{j},\rho_{i+1}]\cdots[\rho_{j},\rho_{j-1}]	&	i<j ,\cr
[[\rho_i,\rho_j],\rho_r]=1											&	i<j<r ,\cr
[[\rho_i,\rho_r^{-1}],\rho_j]=1										&	i<j<r 	 ,\cr
[[\rho_j,\rho_i],[\rho_r,\rho_i][\rho_r,\rho_j]]=1 					& 	i<j<r	,\cr
[[\rho_r,\rho_i],[\rho_r,\rho_j][\rho_j,\rho_i]]=1 					& 	i<j<r,\cr
[[\rho_s,\rho_j],[\rho_j,\rho_r][\rho_r,\rho_i][\rho_r,\rho_j]]=1 	& 	i<j<r<s,\cr
[[\rho_s,\rho_j],[\rho_s,\rho_r][\rho_r,\rho_i][\rho_r,\rho_s]]=1 	& 	i<j<r<s.\cr
\end{array}
\]
\end{itemize}
Moreover, the image of the generator $A_{j,i}\in P_k(A)$, $j>i$, by the map induced by the natural inclusion $A\hookrightarrow \Cal M$ is $[\rho_j,\rho_i]$.
\end{theo}

\begin{demo}
The group $P_k(\Cal M)$ is the fundamental group of the geometrical realization of the following $h\cube_2$-space.
$$\xymatrix{
{\left({\Cal C}_2^k\times{\FF_2(C)}\times \FF_{k-2}(A)\right)^{(0)}} \ar@<-0.5ex>[r]\ar@<0.5ex>[r]\ar@<-1.5ex>[r]\ar@<1.5ex>[r]	&
{\Cal C}_1^k\times C \times (\FF_{k-1}(A))\ar@<-0.5ex>[r]\ar@<0.5ex>[r]&
{\FF_k(A)}.\\
}$$
Let $s=(0,-1)\in\FF_2(C)$ and $t=(0,1)\in\FF_2(C)$.
Up to an homotopy equivalence, the preceding diagram restricts to the following one named $Y_\bullet$.
$$\xymatrix{
{\Cal C}_2^k\times{\{s,t\}}\times \{(q_1,\dots,q_{k-2})\} \ar@<-0.5ex>[r]\ar@<0.5ex>[r]\ar@<-1.5ex>[r]\ar@<1.5ex>[r]	&
{\Cal C}_1^k \times \{0\}\times\FF_{k-1}(A)\ar@<-0.5ex>[r]\ar@<0.5ex>[r]&
{\FF_k(A)}\\
}.$$
\noindent\textbf{First step: }
We want to apply the Van Kampen Theorem \ref{VK-tore} to obtain the fundamental group of $|Y_{\leq 1}|$. 
So we fix $r\in \ito{k}$ and consider the mapping torus
$\xymatrix{
\{r\}\times \{0\}\times\FF_{k-1}(A)\ar@<-0.5ex>[r]\ar@<0.5ex>[r]&
{\FF_k(A)}\\
}$.
Define some paths in $\FF_k(A)$, denoted by $y_r^\epsilon$, moving the base point
$Q_k=(q_1,\dots,q_{r-1},q_r,q_{r+1},\dots,q_{k})$ to the configuration $(q_1,\dots,q_{r-1},i^\epsilon(0),q_{r},\dots,q_{k-1})$ at constant speed along a segment like below.
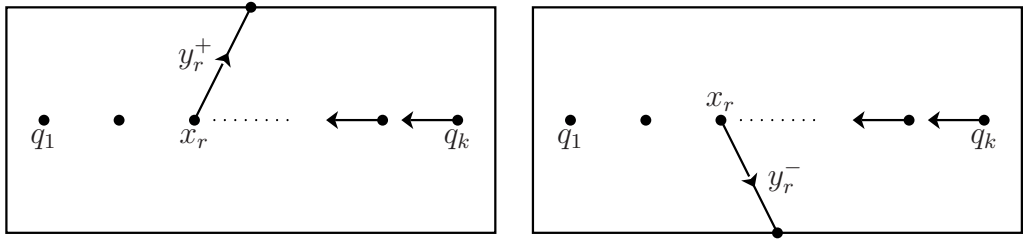
\begin{center}
\psset{xunit=1mm,yunit=1mm,runit=1mm}
\psset{linewidth=0.3,dotsep=1,hatchwidth=0.3,hatchsep=1.5,shadowsize=1}
\psset{dotsize=0.7 2.5,dotscale=1 1,fillcolor=black}
\psset{arrowsize=1 2,arrowlength=1,arrowinset=0.25,tbarsize=0.7 5,bracketlength=0.15,rbracketlength=0.15}
\begin{pspicture}(0,0)(135,35)
\psdots[](15,20)
(5,20)
\psdots[](25,20)
(32.5,35)
\psdots[](50,20)
(60,20)
\pspolygon[](0,5)(65,5)(65,35)(0,35)
\rput[t](5,18.75){$q_1$}
\rput[t](25,18.75){$x_r$}
\rput[t](60,18.75){$q_k$}
\rput[r](27.5,28.75){$y_r^+$}
\psline[arrowscale=1.5 1]{-<}(32.5,35)(28.75,27.5)
\psline(25,20)(28.75,27.5)
\psline[arrowscale=1.5 1]{->}(50,20)(42.5,20)
\psline[arrowscale=1.5 1]{->}(60,20)(52.5,20)
\psline[linestyle=dotted](27.5,20)(37.5,20)
\psdots[](85,20)
(75,20)
\psdots[](95,20)
(102.5,5)
\psdots[](120,20)
(130,20)
\pspolygon[](70,5)(135,5)(135,35)(70,35)
\rput[t](75,18.75){$q_1$}
\rput[b](95,21.25){$x_r$}
\rput[t](130,18.75){$q_k$}
\rput[l](101.25,12.5){$y_r^-$}
\psline[arrowscale=1.5 1]{-<}(102.5,5)(98.75,12.5)
\psline(95,20)(98.75,12.5)
\psline[arrowscale=1.5 1]{->}(120,20)(112.5,20)
\psline[arrowscale=1.5 1]{->}(130,20)(122.5,20)
\psline[linestyle=dotted](97.5,20)(107.5,20)
\end{pspicture}

\end{center}
The image of $A_{j,i}\in \pi_1(\{r\}\times \FF_{k-1}(A))\equiv P_{k-1}$, $1\leq i < j \leq k$, in $\pi_1(\FF_k(A), Q_k)$ by the map $y_r^- * d_1^-(\bullet)*({y_r^-})^{-1}$ is:
\begin{array}[c]\{{lcl}.
A_{j,i}				&\text{if}	&	i<j<r,	\\
A_{j+1,i+1}			&\text{if}	&	r\leq i<j,\\
A_{j+1,i}&\text{if}	&	i<r\leq j.\\
\end{array}\\
We prove the last line just above, the two others being similar.
The loop $y_r^- * d_1^-(A_{j,i}) *({y_r^-})^{-1}$ is given by the following picture
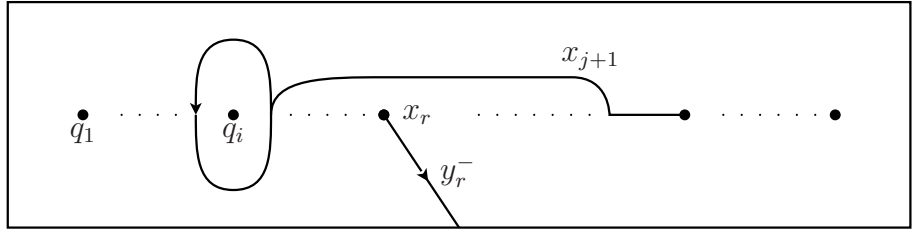
\begin{center}
\psset{xunit=1mm,yunit=1mm,runit=1mm}
\psset{linewidth=0.3,dotsep=1,hatchwidth=0.3,hatchsep=1.5,shadowsize=1}
\psset{dotsize=0.7 2.5,dotscale=1 1,fillcolor=black}
\psset{arrowsize=1 2,arrowlength=1,arrowinset=0.25,tbarsize=0.7 5,bracketlength=0.15,rbracketlength=0.15}
\begin{pspicture}(0,0)(120,30)
\psdots[](10,15)
(30,15)
\psdots[](50,15)
(50,15)
\psdots[](90,15)
(110,15)
\pspolygon[](0,0)(120,0)(120,30)(0,30)
\rput(10,12.5){$q_1$}
\rput[l](52.5,15){$x_r$}
\rput(77.5,22.5){$x_{j+1}$}
\rput[l](57.5,7.5){$y_r^-$}
\psline{-<}(60,0)(55,7.5)
\pscustom[]{\psbezier(90,15)(80,15)(80,15)(80,15)
\psbezier(80,15)(80,15)(80,20)(75,20)
\psbezier(75,20)(70,20)(57.5,20)(50,20)
\psbezier(50,20)(42.5,20)(35,20)(35,15)
}
\rput(30,12.5){$q_i$}
\psline(50,15)(55,7.5)
\psline[linestyle=dotted,dotsep=2](22.5,15)(15,15)
\psline[linestyle=dotted,dotsep=2](47.5,15)(37.5,15)
\psline[linestyle=dotted,dotsep=2](107.5,15)(95,15)
\psline[linestyle=dotted,dotsep=2](77.5,15)(62.5,15)
\pscustom[]{\psbezier{-}(35,15)(35,20)(35,25)(30,25)
\psbezier{->}(30,25)(25,25)(25,20)(25,15)
}
\pscustom[]{\psbezier(25,15)(25,10)(25,5)(30,5)
\psbezier(30,5)(35,5)(35,10)(35,15)
}
\end{pspicture}

\end{center}
and Lemma \ref{commut-lacets} shows that this loop is homotopic to 
$A_{j+1,i}\in P_k$. With a slight adaptation of the proof, we see that the image of $A_{j,i}\in \pi_1(\{r\}\times \FF_{k-1}(A))\equiv P_{k-1}$ in $\pi_1(\FF_k(A), Q_k)$ by the map
$y_r^+ * d_1^+(\bullet)*({y_r^+})^{-1}$ is:
\begin{array}[c]\{{lcl}.
A_{j,i}				&\text{if}	&	i<j<r,	\\
A_{j+1,i+1}			&\text{if}	&	r\leq i<j,\\
A_{j+1,r}^{-1}A_{j+1,i}A_{j+1,r}&\text{if}	&	i<r\leq j.\\
\end{array}\\

Let $z_r$ be the path in $|Y_{\leq 1}|$ with support the segment
$${r}\times \{0\}\times (q_1,\dots,q_{k-1})\times [-1,1]\subset {\Cal C}_1^k\times C \times (\FF_{k-1}(A))\times [-1,1]\subset |Y_{\leq 1}|.$$
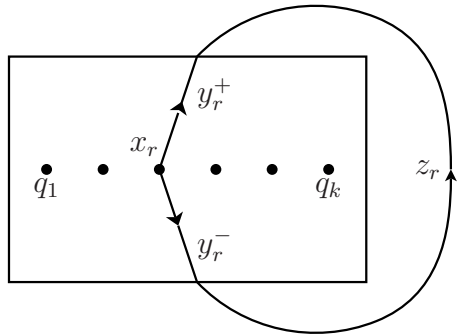
\begin{center}
\psset{xunit=1mm,yunit=1mm,runit=1mm}
\psset{linewidth=0.3,dotsep=1,hatchwidth=0.3,hatchsep=1.5,shadowsize=1}
\psset{dotsize=0.7 2.5,dotscale=1 1,fillcolor=black}
\psset{arrowsize=1 2,arrowlength=1,arrowinset=0.25,tbarsize=0.7 5,bracketlength=0.15,rbracketlength=0.15}
\begin{pspicture}(0,0)(58.75,48.5)
\psdots[](5,25)
(12.5,25)
\psdots[](20,25)
(27.5,25)
\psdots[](35,25)
(42.5,25)
\pspolygon[](0,10)(47.5,10)(47.5,40)(0,40)
\rput(5,22.5){$q_1$}
\rput[Br](20,27.5){$x_r$}
\rput(42.5,22.5){$q_k$}
\pscustom[]{\psbezier(58.75,25)(58.75,33.9)(57.09,41.89)(48.83,45.2)
\psbezier(48.83,45.2)(40.57,48.5)(31.13,46.44)(25,40)
}
\rput[r](57.5,25){$z_r$}
\rput[l](25,15){$y_r^-$}
\rput[l](25,35){$y_r^+$}
\psline[arrowscale=1.5 1]{-<}(25,40)(22.5,32.5)
\psline(20,25)(22.5,32.5)
\psline[arrowscale=1.5 1]{->}(20,25)(22.5,17.5)
\psline(25,10)(22.5,17.5)
\pscustom[]{\psbezier{-}(25,10)(31.13,3.56)(40.57,1.5)(48.83,4.8)
\psbezier{->}(48.83,4.8)(57.09,8.11)(58.75,16.1)(58.75,25)
}
\end{pspicture}

\end{center}
We finally define the loop $\rho_r=y_r^- z_r (y_r^+)^{-1}\in\pi_1(|Y_{\leq 1}|)$.
Proposition \ref{VK-tore} asserts that the space $|Y_{\leq 1}|$
admits the following presentation:
$$\Cal F (\rho_1,\dots,\rho_k)\star P_k / \Cal{R}_1$$
where the relations $\Cal{R}_1$ are given by:
$$
\begin{array}{lr}
		[A_{j,i},\rho_r]=1	&	\textrm{ if } i<j<r	\textrm{ or } r<i<j ,\\
		A_{j,r}^{-1}A_{j,i}A_{j,r}\rho_r^{-1}A_{j,i}^{-1}\rho_r = 1	& 	\textrm{ if } i<r<j .\\
\end{array}
$$

\noindent\textbf{Second step: }
The space $|Y_{\leq2}|$ is obtained by gluing a $2$-cell for each element of $\Cal C_2^k\times\{s,t\}\equiv\Cal C_2^k\times\{s,t\}\times \{(q_1,\dots,q_{k-2})\}$ to the space $|Y_{\leq 1}|$.
Fix $(p,q)\in\Cal C_2^k$. We carry out the details for the $2$-cell associated to the point $((p,q),s)$.\\
First, we recall the paths $\phi_1$, $\phi_2$, $\phi_3$ and $\phi_4$ defined in Example \ref{realization-cube2}.
They are briefly schematized in the following way:
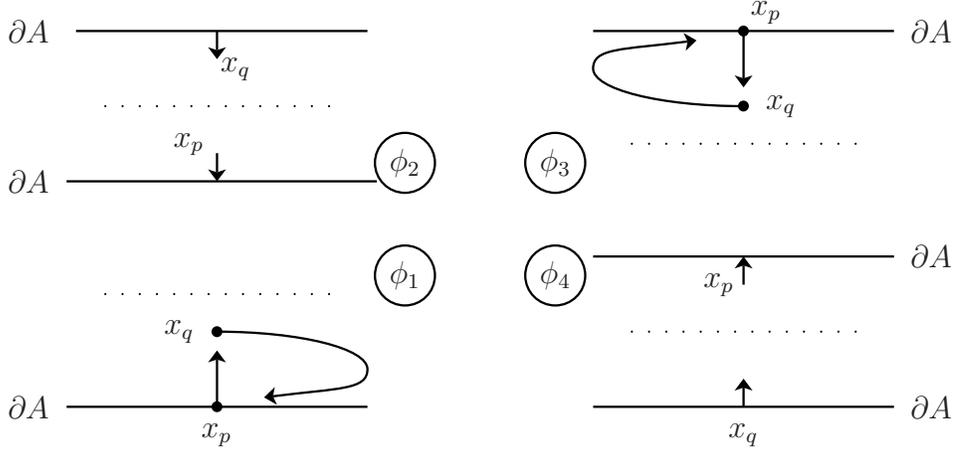
\begin{center}
\psset{xunit=1mm,yunit=1mm,runit=1mm}
\psset{linewidth=0.3,dotsep=1,hatchwidth=0.3,hatchsep=1.5,shadowsize=1}
\psset{dotsize=0.7 2.5,dotscale=1 1,fillcolor=black}
\psset{arrowsize=1 2,arrowlength=1,arrowinset=0.25,tbarsize=0.7 5,bracketlength=0.15,rbracketlength=0.15}
\begin{pspicture}(0,0)(115,60)
\psdots[](20,17.5)
(20,7.5)
\psline(0,7.5)(40,7.5)
\rput(20,17.5){}
\pscustom[arrowscale=1.5 1]{\psbezier{<-}(26.25,8.75)(36.25,10)(40,9.74)(40,12.5)
\psbezier{-}(40,12.5)(40,15.26)(31.05,17.5)(20,17.5)
}
\rput[l](27.5,15){}
\rput(15,17.5){$x_q$}
\psline[arrowscale=1.5 1]{->}(20,7.5)(20,15)
\rput(45,25){\pscirclebox[]{$\phi_1$}}
\psdots[](90,47.5)
(90,57.5)
\psline(110,57.5)(70,57.5)
\rput(90,47.5){}
\pscustom[arrowscale=1.5 1]{\psbezier{<-}(83.75,56.25)(73.75,55)(70,55.26)(70,52.5)
\psbezier{-}(70,52.5)(70,49.74)(78.95,47.5)(90,47.5)
}
\rput[l](82.5,50){}
\rput(95,47.5){$x_q$}
\rput[r](95,60){$x_p$}
\psline[arrowscale=1.5 1]{->}(90,57.5)(90,50)
\rput(65,40){\pscirclebox[]{$\phi_3$}}
\rput(-5,7.5){$\partial A$}
\rput(115,57.5){$\partial A$}
\psline(1.25,57.5)(40,57.5)
\psline(0,37.5)(41.25,37.5)
\psline(70,27.5)(110,27.5)
\psline(110,7.5)(70,7.5)
\psline[arrowscale=1.5 1]{->}(20,57.5)(20,53.75)
\psline[arrowscale=1.5 1]{->}(20,41.25)(20,37.5)
\rput(22.5,52.5){$x_q$}
\rput(16.25,42.5){$x_p$}
\psline[arrowscale=1.5 1]{->}(90,23.75)(90,27.5)
\psline[arrowscale=1.5 1]{->}(90,7.5)(90,11.25)
\rput[r](88.75,23.75){$x_p$}
\rput(90,3.75){$x_q$}
\rput(65,25){\pscirclebox[]{$\phi_4$}}
\rput(45,40){\pscirclebox[]{$\phi_2$}}
\rput(-5,57.5){$\partial A$}
\rput(-5,37.5){$\partial A$}
\rput(115,27.5){$\partial A$}
\rput(115,7.5){$\partial A$}
\rput(20,3.75){$x_p$}
\psline[linestyle=dotted,dotsep=2](5,22.5)(35,22.5)
\psline[linestyle=dotted,dotsep=2](5,47.5)(35,47.5)
\psline[linestyle=dotted,dotsep=2](75,17.5)(105,17.5)
\psline[linestyle=dotted,dotsep=2](75,42.5)(105,42.5)
\end{pspicture}

\end{center}
For the sake of clarity, the paths $\phi_1$ and $\phi_3$ are only pictured between the times $0$ and $\frac{3}{4}$. Also, the points $x_i$ for $i\notin \{p,q\}$, which stay motionless, are represented by a dotted line. Recall that the $2$-cell is attached to the space $|Y_{\leq 1}|$ along the loop
\[  \Phi_2(((p,q),s),\bullet) = \phi_1  D_1^-  \phi_2  D_2^+   \phi_3  D_1^+  \phi_4  D_2^-\]
where the paths $D_*^\epsilon$ are defined like in Example \ref{realization-cube2}.
They are represented in the following way:
\[\begin{array}{cc}
D_2^+=d_2^+((p,q),s)\times id, & D_1^+=d_1^+((p,q),s)\times (-id),\\
D_1^-=d_1^-((p,q),s)\times id, & D_2^-=d_2^-((p,q),s)\times (-id).\\
\end{array}\]
A convenient way to study the homotopy class of this $2$-cell attachment  is to split the composition of path into a composition of well pointed loops.
We define eight paths $\alpha_1,\dots,\alpha_8$ in $\FF_k(A)=Y_0$ with origin the base point $(q_1,\dots,q_k)$ such that the following loop, named $\omega$, is well defined:
\[  (\alpha_1 \phi_1 \alpha_2^{-1}) (\alpha_2 D_1^- \alpha_3^{-1}) (\alpha_3 \phi_2 \alpha_4^{-1})  (\alpha_4 D_2^+ \alpha_5^{-1}) (\alpha_5 \phi_3 \alpha_6^{-1})  (\alpha_6 D_1^+ \alpha_7^{-1})  (\alpha_7 \phi_4  \alpha_8^{-1}) (\alpha_8 D_2^- \alpha_1^{-1})\]
In the rest of the proof, the paths $\alpha_*$ consist in moving the points of the configuration linearly and with constant speed as in the following pictures.
\begin{center}
\psset{xunit=1mm,yunit=1mm,runit=1mm}
\psset{linewidth=0.3,dotsep=1,hatchwidth=0.3,hatchsep=1.5,shadowsize=1}
\psset{dotsize=0.7 2.5,dotscale=1 1,fillcolor=black}
\psset{arrowsize=1 2,arrowlength=1,arrowinset=0.25,tbarsize=0.7 5,bracketlength=0.15,rbracketlength=0.15}
\begin{pspicture}(0,0)(150,65)
\psdots[](42.5,15)
(47.5,15)
\psdots[](52.5,15)
(57.5,15)
\psdots[](62.5,15)
(67.5,15)
\pspolygon[](40,0)(70,0)(70,30)(40,30)
\rput[r](47.5,7.5){$x_p$}
\rput[l](62.5,10){$x_q$}
\psline[arrowscale=1.5 1.5]{->}(47.5,15)(55,0)
\psline[arrowscale=1.5 1.5]{->}(62.5,15)(55,7.5)
\rput(65,25){\psframebox[]{$\alpha_1$}}
\psdots[](2.5,15)
(7.5,15)
\psdots[](12.5,15)
(17.5,15)
\psdots[](22.5,15)
(27.5,15)
\pspolygon[](0,0)(30,0)(30,30)(0,30)
\psline[arrowscale=1.5 1.5]{->}(7.5,15)(15,7.5)
\psline[arrowscale=1.5 1.5]{->}(22.5,15)(15,0)
\rput(25,25){\psframebox[]{$\alpha_2$}}
\psdots[](2.5,50)
(7.5,50)
\psdots[](12.5,50)
(17.5,50)
\psdots[](22.5,50)
(27.5,50)
\pspolygon[](0,35)(30,35)(30,65)(0,65)
\psline[arrowscale=1.5 1.5]{->}(7.5,50)(15,42.5)
\psline[arrowscale=1.5 1.5]{->}(22.5,50)(15,65)
\rput(25,40){\psframebox[]{$\alpha_3$}}
\psdots[](42.5,50)
(47.5,50)
\psdots[](52.5,50)
(57.5,50)
\psdots[](62.5,50)
(67.5,50)
\pspolygon[](40,35)(70,35)(70,65)(40,65)
\psline[arrowscale=1.5 1.5]{->}(47.5,50)(55,35)
\psline[arrowscale=1.5 1.5]{->}(62.5,50)(55,57.5)
\rput(65,40){\psframebox[]{$\alpha_4$}}
\psdots[](82.5,15)
(87.5,15)
\psdots[](92.5,15)
(97.5,15)
\psdots[](102.5,15)
(107.5,15)
\pspolygon[](80,0)(110,0)(110,30)(80,30)
\psline[arrowscale=1.5 1.5]{->}(87.5,15)(95,30)
\psline[arrowscale=1.5 1.5]{->}(102.5,15)(95,7.5)
\rput(105,25){\psframebox[]{$\alpha_8$}}
\psdots[](82.5,50)
(87.5,50)
\psdots[](92.5,50)
(97.5,50)
\psdots[](102.5,50)
(107.5,50)
\pspolygon[](80,35)(110,35)(110,65)(80,65)
\psline[arrowscale=1.5 1.5]{->}(87.5,50)(95,65)
\psline[arrowscale=1.5 1.5]{->}(102.5,50)(95,57.5)
\rput(105,40){\psframebox[]{$\alpha_5$}}
\psdots[](122.5,50)
(127.5,50)
\psdots[](132.5,50)
(137.5,50)
\psdots[](142.5,50)
(147.5,50)
\pspolygon[](120,35)(150,35)(150,65)(120,65)
\psline[arrowscale=1.5 1.5]{->}(127.5,50)(135,57.5)
\psline[arrowscale=1.5 1.5]{->}(142.5,50)(135,65)
\rput(145,40){\psframebox[]{$\alpha_6$}}
\psdots[](122.5,15)
(127.5,15)
\psdots[](132.5,15)
(137.5,15)
\psdots[](142.5,15)
(147.5,15)
\pspolygon[](120,0)(150,0)(150,30)(120,30)
\psline[arrowscale=1.5 1.5]{->}(127.5,15)(135,22.5)
\psline[arrowscale=1.5 1.5]{->}(142.5,15)(135,0)
\rput(145,25){\psframebox[]{$\alpha_7$}}
\rput[r](7.5,7.5){$x_p$}
\rput[l](22.5,7.5){$x_q$}
\rput[r](7.5,45){$x_p$}
\rput[l](22.5,57.5){$x_q$}
\rput[r](47.5,42.5){$x_p$}
\rput[l](62.5,55){$x_q$}
\rput[r](87.5,57.5){$x_p$}
\rput[l](102.5,55){$x_q$}
\rput[r](127.5,55){$x_p$}
\rput[l](142.5,57.5){$x_q$}
\rput[r](127.5,20){$x_p$}
\rput[l](142.5,7.5){$x_q$}
\rput[r](87.5,22.5){$x_p$}
\rput[l](102.5,10){$x_q$}
\end{pspicture}
\\
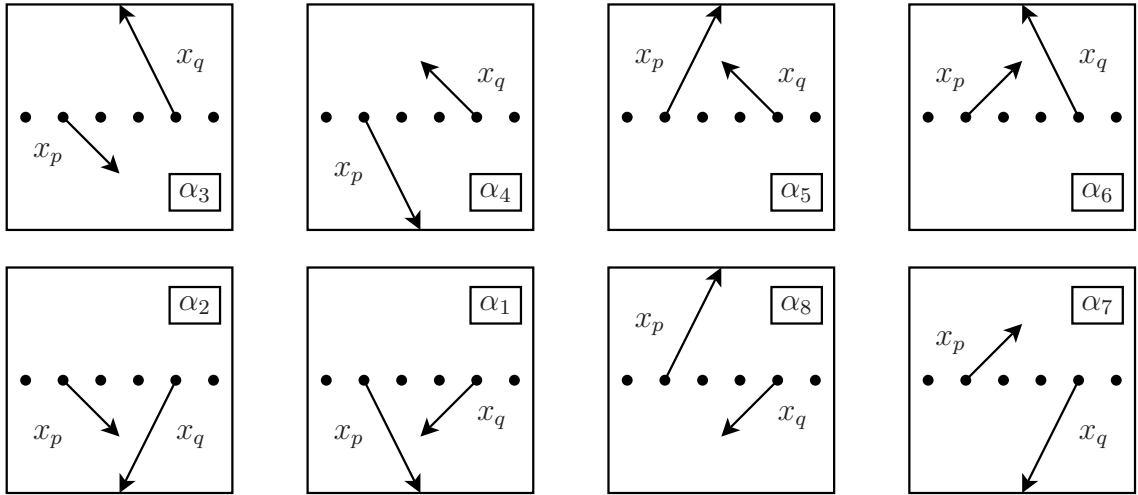
	\small The $8$ paths used for splitting $\omega$
\end{center}
The final step for the determination of the homotopy class $\omega$ consists in drawing the various loops arising from the previous splitting.
For example, consider the loop $\alpha_5 \phi_3 \alpha_6^{-1}$:
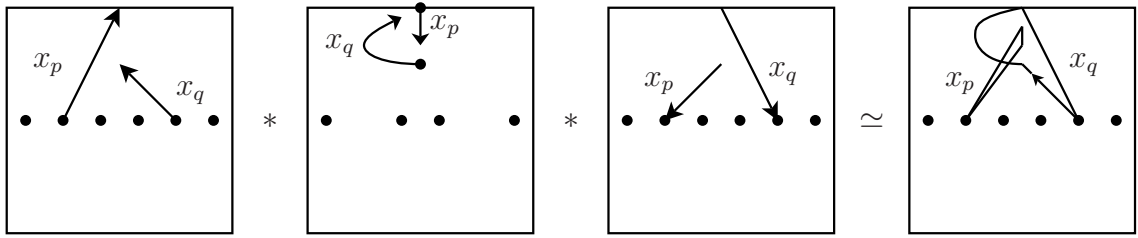
\begin{center}
\psset{xunit=1mm,yunit=1mm,runit=1mm}
\psset{linewidth=0.3,dotsep=1,hatchwidth=0.3,hatchsep=1.5,shadowsize=1}
\psset{dotsize=0.7 2.5,dotscale=1 1,fillcolor=black}
\psset{arrowsize=1 2,arrowlength=1,arrowinset=0.25,tbarsize=0.7 5,bracketlength=0.15,rbracketlength=0.15}
\begin{pspicture}(0,0)(150,30)
\psdots[](122.5,15)
(127.5,15)
\psdots[](132.5,15)
(137.5,15)
\psdots[](142.5,15)
(147.5,15)
\pspolygon[](120,0)(150,0)(150,30)(120,30)
\rput[r](128.75,20){$x_p$}
\rput[l](141.25,22.5){$x_q$}
\psline(127.5,15)(135,25)
\psline{->}(142.5,15)(136.25,21.25)
\psline(127.5,15)(135,27.5)
\psline(142.5,15)(135,30)
\psdots[](2.5,15)
(7.5,15)
\psdots[](12.5,15)
(17.5,15)
\psdots[](22.5,15)
(27.5,15)
\pspolygon[](0,0)(30,0)(30,30)(0,30)
\rput[r](7.5,22.5){$x_p$}
\rput[l](22.5,18.75){$x_q$}
\psline[arrowscale=1.5 1.5]{->}(7.5,15)(15,30)
\psline[arrowscale=1.5 1.5]{->}(22.5,15)(15,22.5)
\psdots[](52.5,15)
(57.5,15)
\psdots[](42.5,15)
(67.5,15)
\pspolygon[](40,0)(70,0)(70,30)(40,30)
\rput[l](56.25,27.5){$x_p$}
\rput[r](46.25,25){$x_q$}
\psdots[](82.5,15)
(87.5,15)
\psdots[](92.5,15)
(97.5,15)
\psdots[](102.5,15)
(107.5,15)
\pspolygon[](80,0)(110,0)(110,30)(80,30)
\rput[r](88.75,20){$x_p$}
\rput[l](101.25,21.25){$x_q$}
\psline[arrowscale=1.5 1.5]{->}(95,22.5)(87.5,15)
\psline[arrowscale=1.5 1.5]{->}(95,30)(102.5,15)
\psdots[](55,30)
(55,22.5)
\rput(56.25,5){}
\pscustom[arrowscale=1.5 1]{\psbezier{<-}(52.5,28.75)(50,27.5)(47.5,26.25)(47.5,25)
\psbezier{-}(47.5,25)(47.5,23.75)(50,22.5)(55,22.5)
}
\rput[l](60,3.75){}
\psline[arrowscale=1.5 1]{->}(55,30)(55,25)
\rput(121.25,-1.25){}
\rput[l](125,-2.5){}
\psline(135,25)(135,27.5)
\rput(35,15){$\ast$}
\rput(75,15){$\ast$}
\rput(115,15){$\he$}
\pscustom[]{\psbezier(135,30)(128.75,28.75)(128.75,27.5)(128.75,26.25)
\psbezier(128.75,26.25)(128.75,25)(130,22.5)(135,22.5)
}
\psline(135,22.5)(136.25,21.25)
\end{pspicture}

\end{center}
To find its homotopy class, observe that $\alpha_5$ and $\alpha_6^{-1}$ can be re-parametrized such that the $p$-th point is moved before the $q$-th one.
Hence, $\alpha_5 \phi_3 \alpha_6^{-1}\he A_{q,p}^{-1}$ as pictured below.
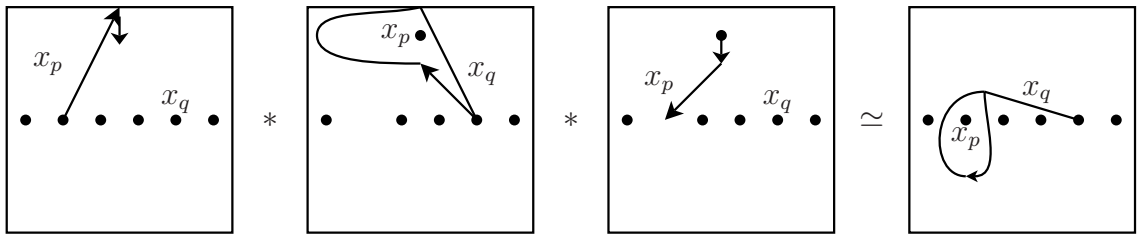
\begin{center}
\psset{xunit=1mm,yunit=1mm,runit=1mm}
\psset{linewidth=0.3,dotsep=1,hatchwidth=0.3,hatchsep=1.5,shadowsize=1}
\psset{dotsize=0.7 2.5,dotscale=1 1,fillcolor=black}
\psset{arrowsize=1 2,arrowlength=1,arrowinset=0.25,tbarsize=0.7 5,bracketlength=0.15,rbracketlength=0.15}
\begin{pspicture}(0,0)(150,30)
\psdots[](122.5,15)
(127.5,15)
\psdots[](132.5,15)
(137.5,15)
\psdots[](142.5,15)
(147.5,15)
\pspolygon[](120,0)(150,0)(150,30)(120,30)
\rput(127.5,12.5){$x_p$}
\rput[l](135,18.75){$x_q$}
\psline(142.5,15)(130,18.75)
\psdots[](2.5,15)
(7.5,15)
\psdots[](12.5,15)
(17.5,15)
\psdots[](22.5,15)
(27.5,15)
\pspolygon[](0,0)(30,0)(30,30)(0,30)
\rput[r](7.5,22.5){$x_p$}
\rput(22.5,17.5){$x_q$}
\psline[arrowscale=1.5 1.5]{->}(7.5,15)(15,30)
\psline[arrowscale=1.5 1.5]{->}(62.5,15)(55,22.5)
\psdots[](52.5,15)
(57.5,15)
\psdots[](42.5,15)
(67.5,15)
\pspolygon[](40,0)(70,0)(70,30)(40,30)
\rput[r](53.75,26.25){$x_p$}
\rput[l](61.25,21.25){$x_q$}
\psdots[](82.5,15)
(95,26.25)
\psdots[](92.5,15)
(97.5,15)
\psdots[](102.5,15)
(107.5,15)
\pspolygon[](80,0)(110,0)(110,30)(80,30)
\rput[r](88.75,20){$x_p$}
\rput(102.5,17.5){$x_q$}
\psline[arrowscale=1.5 1.5]{->}(95,22.5)(87.5,15)
\psline(55,30)(62.5,15)
\psdots[](55,26.25)
(62.5,15)
\rput(56.25,5){}
\pscustom[]{\psbezier(55,30)(50,28.75)(41.25,30)(41.25,26.25)
\psbezier(41.25,26.25)(41.25,22.5)(50,22.5)(55,22.5)
}
\rput[l](60,3.75){}
\psline[arrowscale=1.5 1]{->}(15,28.75)(15,25)
\rput(121.25,-1.25){}
\rput[l](125,-2.5){}
\rput(35,15){$\ast$}
\rput(75,15){$\ast$}
\rput(115,15){$\he$}
\psbezier(130,18.75)(122.5,18.75)(122.5,7.5)(127.5,7.5)
\psline[arrowscale=1.5 1]{->}(95,26.25)(95,22.5)
\psbezier{<-}(127.5,7.5)(132.5,7.5)(130,15)(130,18.75)
\end{pspicture}

\end{center}
A direct application of  Lemma \ref{commut-lacets} shows that the three loops $\alpha_1 \phi_1 \alpha_2^{-1}$, $\alpha_3 \phi_2 \alpha_4^{-1}$ and $\alpha_7 \phi_4  \alpha_8^{-1}$ are all nullhomotopic.

Also the maps $d_1^-$ sends $((p,q),s)$ in $\{q\}\times\FF_{k-1}(A)$.
The loop $\alpha_2 D_1^- \alpha_3^{-1}$ is pictured below.
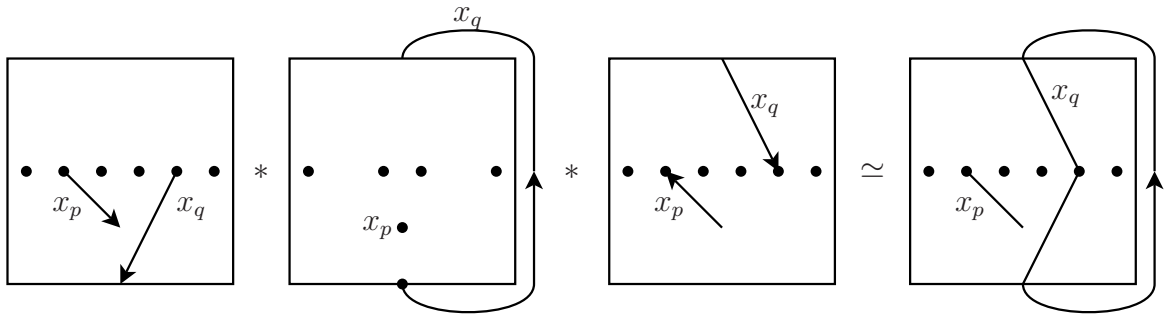
\begin{center}
\psset{xunit=1mm,yunit=1mm,runit=1mm}
\psset{linewidth=0.3,dotsep=1,hatchwidth=0.3,hatchsep=1.5,shadowsize=1}
\psset{dotsize=0.7 2.5,dotscale=1 1,fillcolor=black}
\psset{arrowsize=1 2,arrowlength=1,arrowinset=0.25,tbarsize=0.7 5,bracketlength=0.15,rbracketlength=0.15}
\begin{pspicture}(0,0)(152.5,42.5)
\psdots[](122.5,20)
(127.5,20)
\psdots[](132.5,20)
(137.5,20)
\psdots[](142.5,20)
(147.5,20)
\pspolygon[](120,5)(150,5)(150,35)(120,35)
\rput[r](130,15){$x_p$}
\rput[l](138.75,30){$x_q$}
\psline(127.5,20)(135,12.5)
\psline(142.5,20)(135,5)
\psdots[](2.5,20)
(7.5,20)
\psdots[](12.5,20)
(17.5,20)
\psdots[](22.5,20)
(27.5,20)
\pspolygon[](0,5)(30,5)(30,35)(0,35)
\rput[r](10,15){$x_p$}
\rput[l](22.5,15){$x_q$}
\psline[arrowscale=1.5 1.5]{->}(7.5,20)(15,12.5)
\psline[arrowscale=1.5 1.5]{->}(22.5,20)(15,5)
\psdots[](50,20)
(55,20)
\psdots[](40,20)
(65,20)
\pspolygon[](37.5,5)(67.5,5)(67.5,35)(37.5,35)
\rput[r](51.25,12.5){$x_p$}
\rput[B](61.25,40){$x_q$}
\psdots[](82.5,20)
(87.5,20)
\psdots[](92.5,20)
(97.5,20)
\psdots[](102.5,20)
(107.5,20)
\pspolygon[](80,5)(110,5)(110,35)(80,35)
\rput[r](90,15){$x_p$}
\rput[l](98.75,28.75){$x_q$}
\psline[arrowscale=1.5 1.5]{->}(95,12.5)(87.5,20)
\psline[arrowscale=1.5 1.5]{->}(95,35)(102.5,20)
\psdots[](52.5,12.5)
(52.5,5)
\rput(53.75,10){}
\rput[l](57.5,8.75){}
\rput(121.25,3.75){}
\rput[l](125,2.5){}
\rput(33.75,20){$\ast$}
\rput(75,20){$\ast$}
\rput(115,20){$\he$}
\pscustom[arrowscale=1.5 1.5]{\psbezier{-}(52.5,5)(52.5,0)(70,0)(70,5)
\psbezier{->}(70,5)(70,10)(70,15)(70,20)
}
\pscustom[]{\psbezier(52.5,35)(52.5,40)(70,40)(70,35)
\psbezier(70,35)(70,30)(70,25)(70,20)
}
\pscustom[arrowscale=1.5 1.5]{\psbezier{-}(135,5)(135,0)(152.5,0)(152.5,5)
\psbezier{->}(152.5,5)(152.5,10)(152.5,15)(152.5,20)
}
\pscustom[]{\psbezier(135,35)(135,40)(152.5,40)(152.5,35)
\psbezier(152.5,35)(152.5,30)(152.5,25)(152.5,20)
}
\psline(135,35)(142.5,20)
\end{pspicture}

\end{center}
Observe that up to the position of the $p$-th particle, $D_1^-$ is the path $z_q$.
Hence, using the lemma \ref{commut-lacets}, we come to the conclusion that
$\alpha_2 D_1^- \alpha_3^{-1}\he \rho_q$.
The cases $(\alpha_4 D_2^+ \alpha_5^{-1})\he \rho_p$, $(\alpha_6 D_1^+ \alpha_7^{-1})\he\rho_q^{-1}$ and $(\alpha_8 D_2^- \alpha_1^{-1})\he \rho_p^{-1}$ are obtained analogously.
Finally, the homotopy class of the loop $\omega$ is
$ \rho_q\rho_p A_{q,p}^{-1}\rho_q^{-1}\rho_p^{-1}$.

In a similar way, the $2$-cell associated to the element $((p,q),t)$, once split and pointed at $Q_k$, is attached along the loop
$B_{q,p}*\rho_q*1*\rho_p*1*\rho_q^{-1}*1*\rho_p^{-1}\he B_{q,p}\rho_q\rho_p\rho_q^{-1}\rho_p^{-1}$.

As proved in Section \ref{Method}, the group $P_k(\Cal M)=\pi_1(|Y_{\leq2}|,Q_k)$ admits the presentation:
\[\Cal F (\rho_1,\dots,\rho_k)\star P_k / \Cal{R}_2\]
where the relations $\Cal{R}_2$ are given by:
\[
\begin{array}[t]{lll}
(a)	&		[A_{j,i},\rho_r]=1	&	\textrm{ if } i<j<r	\textrm{ or } r<i<j, \\
(b)&		A_{j,r}^{-1}A_{j,i}A_{j,r}\rho_r^{-1}A_{j,i}^{-1}\rho_r = 1	& 	\textrm{ if } i<r<j,\\
(c)	&	B_{q,p}=[\rho_p^{-1},\rho_q^{-1}]& 	\textrm{ if } p<q ,\\
	&	A_{q,p}=[\rho_q,\rho_p]& 	\textrm{ if } p<q .\\
\end{array}
\]

The final step of the proof consists of a simplification of the presentation.
For that, we first replace $A_{q,p}$ with $[\rho_q,\rho_p]$ in the previous relations.
\begin{itemize}
    \item[(a)] If $i<j<r$  or   $r<i<j$, then $[A_{j,i},\rho_r]=[[\rho_j,\rho_i],\rho_r]=1$.
	\item[(b)] If $i<r<j$, then $A_{j,r}^{-1}A_{j,i}A_{j,r}\rho_r^{-1}A_{j,i}^{-1}\rho_r =\rho_r^{-1}\rho_j^{-1}(\rho_r\rho_i^{-1}\rho_j\rho_i\rho_j^{-1}\rho_r^{-1}\rho_j\rho_i^{-1}\rho_j^{-1}\rho_i)\rho_j\rho_r =1$\\[.5em] which implies \quad $\rho_r\rho_i^{-1}\rho_j\rho_i\rho_j^{-1}\rho_r^{-1}\rho_j\rho_i^{-1}\rho_j^{-1}\rho_i =1$ \quad i.e. \quad$ [\rho_r,[\rho_i,\rho_j^{-1}]]=1$.
	\item[(c)] The relation
$B_{j,i}=A_{j,j-1}^{-1}A_{j,j-2}^{-1}\cdots A_{j,i+1}^{-1} A_{j,i} A_{j,i+1}\cdots A_{j,j-2} A_{j,j-1}$ \quad 
in $P_k$ implies
$$
[\rho_i^{-1},\rho_j^{-1}]=[\rho_{j-1},\rho_j][\rho_{j-2},\rho_j]\cdots [\rho_{i+1},\rho_j]
[\rho_{j},\rho_i][\rho_{j},\rho_{i+1}]\cdots[\rho_{j},\rho_{j-2}][\rho_{j},\rho_{j-1}].
$$
\end{itemize}
Observe that, as a consequence of relation $(a)$, we have
$[[\rho_j,\rho_i],[\rho_s,\rho_r]]=[A_{j,i},A_{s,r}]=1$ for $i<j<r<s$ or $r<i<j<s$.
In fact, these are two of the Yang-Baxter relations given in \ref{equ-tresses} ($(3)$ and $(4)$).
Finally, the remaining Yang-Baxter relations translate into the following relations.
\[
\begin{array}{rcll}
[A_{j,i},A_{r,i}A_{r,j}]=1	&\Rightarrow& [[\rho_j,\rho_i],[\rho_r,\rho_i][\rho_r,\rho_j]]=1 &\textrm{if }  i<j<r ,	\cr
[A_{r,i},A_{r,j}A_{j,i}]=1	&\Rightarrow&
[[\rho_r,\rho_i],[\rho_r,\rho_j][\rho_j,\rho_i]]=1 &\textrm{if } i<j<r ,\cr
[A_{s,j},A_{r,j}^{-1}A_{r,i}A_{r,j}]=1	&\Rightarrow&
[[\rho_s,\rho_j],[\rho_j,\rho_r][\rho_r,\rho_i][\rho_r,\rho_j]]=1 &\textrm{if } i<j<r<s ,\cr
[A_{s,j},A_{s,r}^{-1}A_{r,i}A_{s,r}]=1	&\Rightarrow&
[[\rho_s,\rho_j],[\rho_s,\rho_r][\rho_r,\rho_i][\rho_r,\rho_s]]=1&\textrm{if } i<j<r<s .\cr
\end{array}
\]
\end{demo}

Of course, we can iterate the construction by attaching a dimension $2$ disc $D^2$ along the boundary of the Möbius band in order to obtain the projective plane $\RR P^2$.
In the particular case $k=2$, we recover a result from J.~van~Buskirk\cite{van1966}.
\begin{CORO}[\cite{van1966}]
The pure braid group on $2$ strands of the projective plane $\RR P^2$ is isomorphic to the group of the quaternions $Q_8$.
\end{CORO}

\small

\bibliographystyle{abbrv}
\bibliography{../../biblio}

\noindent
Université de Lille 1\\ Département de Mathématiques --- UMR 8524\\ 59655 Villeneuve d'Ascq Cedex, France\\
\texttt{jourdan@math.univ-lille1.fr}

\end{document}